\documentclass[11pt, reqno]{amsart}

\usepackage[a4paper,margin=1in]{geometry}
\usepackage[T1]{fontenc}
\usepackage[utf8]{inputenc}
\IfFileExists{lmodern.sty}{\usepackage{lmodern}}{}
\IfFileExists{microtype.sty}{\usepackage{microtype}}{}
\usepackage{amsmath,amssymb,amsthm,mathtools,bm,mathrsfs}
\usepackage{enumitem}
\usepackage{hyperref,xcolor}
\IfFileExists{kotex.sty}{\usepackage{kotex}}{}
\usepackage[normalem]{ulem}

\makeatletter
\renewcommand{\paragraph}{\@startsection{paragraph}{4}%
  {\z@}{.5\baselineskip}{-.5em}{\normalfont\bfseries}}
\makeatother

\allowdisplaybreaks

\numberwithin{equation}{section}

\newtheorem{theorem}{Theorem}[section]
\newtheorem{proposition}[theorem]{Proposition}
\newtheorem{lemma}[theorem]{Lemma}
\newtheorem{corollary}[theorem]{Corollary}
\theoremstyle{definition}
\newtheorem{definition}[theorem]{Definition}
\theoremstyle{remark}
\newtheorem{remark}[theorem]{Remark}

\newcommand{\R}{\mathbb{R}}
\newcommand{\Rd}{\mathbb{R}^d}
\newcommand{\dd}{\,\mathrm{d}}
\newcommand{\dist}{\operatorname{dist}}
\newcommand{\proj}{\operatorname{proj}}
\newcommand{\mcP}{\mathcal{P}}
\newcommand{\one}{\mathbf{1}}
\newcommand{\eps}{\varepsilon}
\newcommand{\Tr}{\operatorname{Tr}}
\newcommand{\divg}{\nabla\!\cdot\!}
\newcommand{\divG}{\operatorname{div}_{\Gamma}}
\newcommand{\nG}{\nabla_{\!\Gamma}}
\newcommand{\wl}{w_\lambda}
\newcommand{\muh}{\widehat\mu_\lambda}
\newcommand{\hh}{\widehat h}
\newcommand{\Hm}{\mathcal{H}^{d-1}}
\newcommand{\dBL}{d_{\mathrm{BL}}}
\newcommand{\norm}[1]{\left\lVert #1\right\rVert}
\newcommand{\abs}[1]{\left\lvert #1\right\rvert}

\title{Soft-confinement approximation of no-flux McKean--Vlasov dynamics}
\author[Geuntaek Seo]{Geuntaek Seo\textsuperscript{*}}
\address{Department of Mathematics, Pohang University of Science and Technology}
\email{gtseo@postech.ac.kr}
\thanks{\textsuperscript{*}Geuntaek Seo and Yeoneung Kim contributed equally to this work.}

\author[Yeoneung Kim]{Yeoneung Kim\textsuperscript{*}}
\address{Department of Industrial Engineering, Yonsei University}
\email{yeoneung@yonsei.ac.kr}

\author{Beomjun Choi}
\address{Department of Mathematical Sciences, Korea Advanced Institute of Science and Technology}
\email{bchoi@kaist.ac.kr}
\hypersetup{pdfauthor={Geuntaek Seo, Yeoneung Kim, and Beomjun Choi}}

\begin{document}

\begin{abstract}
We study the approximation of the no-flux McKean--Vlasov equation on a bounded convex domain by a whole-space equation with quadratic confinement of strength $\lambda^{-1}$ outside the domain. For a smooth bounded interaction kernel and a common bounded initial density with square-integrable gradient, we establish an $O(\lambda^{1/4})$ approximation rate. The estimate controls the $L^2$ density error uniformly on every finite time interval and the gradient error in space--time $L^2$. Without interaction, a spectral gap independent of $\lambda$ yields the $O(\lambda^{1/4})$ rate uniformly in time for the density error and uniformly over all unit time intervals for the gradient error. With interaction, we construct a counterexample showing that even weak convergence can fail to be uniform in time.
\end{abstract}



\maketitle

\section{Introduction}\label{sec:intro}

McKean--Vlasov equations describe the competition between diffusion, an external potential, and interactions among particles. This structure appears in synchronization and chemotaxis, including the noisy Kuramoto and parabolic--elliptic Keller--Segel models \cite{CarrilloGvalaniPavliotisSchlichting2020}. On the sphere, related equations describe orientational alignment in Doi--Smoluchowski models \cite{ConstantinKevrekidisTiti2004} and noisy self-attention dynamics \cite{ShalovaSchlichting2026}. In sampling and optimization, entropy-regularized mean-field dynamics combine diffusion with descent of an energy on probability measures \cite{Chizat2022}. The training of wide two-layer neural networks also admits a mean-field description as a Wasserstein gradient flow of the parameter distribution \cite{MeiMontanariNguyen2018,ChizatBach2018}. For squared loss, adding Gaussian noise to the parameter dynamics gives a McKean--Vlasov diffusion \cite{MeiMontanariNguyen2018,HuRenSiskaSzpruch2021}. A notable recent example is the entropy-regularized maximum mean discrepancy (MMD) flow: a kernel determines both interaction among particles and attraction to a target distribution \cite{ArbelKorbaSalimGretton2019,SeoKimChoi2026}.

We study how accurately a strong confining potential approximates a no-flux boundary condition for McKean--Vlasov equations with smooth bounded interaction kernels on a bounded convex Euclidean domain. Reflection keeps particle trajectories in the domain, and the no-flux condition expresses this constraint at the density level \cite{Sznitman1984,HindsSharmaTretyakov2024}. The penalized equation allows mass to leave the domain and drives it back with a restoring force. Our aim is to quantify the difference between these two evolutions as the confinement becomes stronger, and to determine whether the approximation remains valid uniformly in time.

Let \(K\subset\Rd\) be compact and convex, with nonempty interior and \(C^3\) boundary. Equations and Sobolev spaces on \(K\) are understood on its interior. The constrained evolution is
\begin{equation}\label{eq:intro-no-flux}
  \begin{cases}
    \partial_t\rho=\divg\bigl(\nabla\rho+\rho b[\rho(t)]\bigr),
      &\text{in }(0,\infty)\times K,\\
    \bigl(\nabla\rho+\rho b[\rho(t)]\bigr)\cdot\vec n=0,
      &\text{on }(0,\infty)\times\partial K,
  \end{cases}
\end{equation}
where \(\vec n\) is the outward unit normal and
\[
  b[\mu](x)=\nabla\Phi(x)
    +\int_{\Rd}\nabla_xW(x,y)\,\dd\mu(y).
\]
The function \(\Phi\) is the external potential, and the symmetric kernel \(W(x,y)\) is the pairwise interaction potential. Measures on \(K\) are regarded as measures on \(\Rd\) supported in \(K\). The associated free energy is
\begin{equation}\label{eq:intro-free-energy}
  \mathcal F_K(\rho)
  =\int_K\rho\log\rho\,\dd x
   +\int_K\Phi\rho\,\dd x
   +\frac12\iint_{K\times K}W(x,y)\rho(x)\rho(y)\,\dd x\dd y.
\end{equation}
These three terms are the entropy, external potential energy, and interaction energy, respectively. Up to an additive constant, its first variation is
\[
  \frac{\delta\mathcal F_K}{\delta\rho}(x)
  =\log\rho(x)+\Phi(x)+\int_K W(x,y)\rho(y)\,\dd y.
\]
Thus \eqref{eq:intro-no-flux} has the form \(\partial_t\rho=\divg(\rho\nabla(\delta\mathcal F_K/\delta\rho))\), the Wasserstein gradient flow of \(\mathcal F_K\) on probability measures supported in \(K\).  %Formally, the no-flux condition gives
%\[
%  \frac{\dd}{\dd t}\mathcal F_K(\rho(t))
%  =-\int_K\rho\left|\nabla\frac{\delta\mathcal F_K}{\delta\rho}\right|^2\dd x.
%\]

Soft confinement replaces the support constraint by the additional potential
\[
  U_\lambda(x)=\frac{\dist(x,K)^2}{2\lambda},
  \qquad \lambda>0,
\]
the Moreau--Yosida regularization of the indicator of \(K\) \cite{Moreau1965}. Keeping \(\Phi\) and \(W\) fixed, we consider
\begin{equation}\label{eq:intro-penalized}
  \partial_t\rho^\lambda
  =\divg\left(\nabla\rho^\lambda+\rho^\lambda\nabla U_\lambda
    +\rho^\lambda b[\rho^\lambda(t)]\right)
  \qquad\text{in }(0,\infty)\times\Rd.
\end{equation}
Its free energy is obtained from \eqref{eq:intro-free-energy} by taking the integrals over \(\Rd\) and adding \(\int_{\Rd}U_\lambda\rho\,\dd x\). The restoring drift \(-\nabla U_\lambda(x)=-\lambda^{-1}(x-\proj_Kx)\) vanishes on \(K\) and pushes exterior mass toward its nearest point in \(K\). This penalization is also used in constrained Langevin sampling \cite{BrosseDurmusMoulinesPereyra2017,GurbuzbalabanHuZhu2024}. We start both evolutions from a common probability density \(\rho_0\) on \(K\), extended by zero for \eqref{eq:intro-penalized}.

Our first result shows that, as \(\lambda\downarrow0\), the penalized density approximates the no-flux solution with error \(O(\lambda^{1/4})\) in \(L^\infty(0,T;L^2(K))\cap L^2(0,T;H^1(K))\) for every finite \(T\). We obtain the same rate for the interior flux and quantify the decay of its normal component on \(\partial K\). When \(W=0\), a uniform spectral gap gives this rate for all times, with gradient and flux estimates on every unit time interval. With interaction, however, even weak convergence can fail to be uniform in time. We construct an example in which the no-flux solution stays at an unstable equilibrium, while the small amount of mass escaping through the soft boundary creates a perturbation that grows to a size independent of \(\lambda\).

\subsection{Main results}\label{subsec:main-results}

Throughout the paper, \(K\) is as above, \(\Phi\in C^2(\Rd)\) with \(\nabla\Phi\in C_b^1(\Rd;\Rd)\), and \(W\in C_b^2(\Rd\times\Rd)\) is symmetric. Unless otherwise stated, \(\rho_0\in L^\infty(K)\) is a probability density, and the initial data are \(\rho_0\) for the no-flux equation and \(\rho_0\one_K\) for the penalized equation. We call this choice \emph{hard initialization}. The weak formulations are given in Section~\ref{subsec:weak-formulations}.

Let \(\rho^\lambda\) and \(\rho\) denote these solutions, and write their interior fluxes as
\begin{equation}\label{eq:JK-def}
  J_K^\lambda:=\nabla\rho^\lambda+\rho^\lambda b[\rho^\lambda(t)],
  \qquad
  J:=\nabla\rho+\rho b[\rho(t)].
\end{equation}
Norms on \(K\) refer to the restriction of \(\rho^\lambda\) to \(K\).

\begin{theorem}[Strong convergence on bounded time intervals]
\label{thm:strong-main}
Suppose additionally that \(\rho_0\in H^1(K)\). For every \(0<T<\infty\) there are \(C_T>0\) and \(\lambda_0>0\) such that
\begin{equation}\label{eq:strong-rate}
  \begin{aligned}
  &\sup_{0\le t\le T}\|\rho^\lambda(t)-\rho(t)\|_{L^2(K)}
  +\|\rho^\lambda-\rho\|_{L^2(0,T;H^1(K))}+\|J^\lambda_K-J\|_{L^2((0,T)\times K)}
  \le C_T\lambda^{1/4}
  \end{aligned}
\end{equation}
for \(0<\lambda\le\lambda_0\). Moreover, for every \(\tau\in(0,T)\),
\begin{equation}\label{eq:bdry-rate}
  \|J^\lambda_K\cdot\vec n\|_{L^2(\tau,T;L^2(\partial K))}
  \le C_T(1+\tau^{-1/2})\lambda^{1/4}.
\end{equation}
If the hard initialization is replaced by the $\lambda$-dependent auxiliary initial densities constructed in Section~\ref{subsec:lifted}, the corresponding trace estimate holds on the full interval \((0,T)\). More precisely, writing \(\bar J_K^\lambda\) for the interior flux of the auxiliary solution, we have
\[
  \|\bar J_K^\lambda\cdot\vec n\|_{L^2(0,T;L^2(\partial K))}
  \le C_T\lambda^{1/4}.
\]
\end{theorem}
\begin{remark}
The exponent $1/4$ in \eqref{eq:strong-rate} is optimal for each of the three error terms. Matching lower bounds hold already when $K=[0,1]$, $\Phi=W=0$, and $\rho_0\equiv1$; see Proposition~\ref{prop:sharpness}.
\end{remark}

The first estimate controls both the density error at each time and the space--time gradient and flux errors. The second describes how the no-flux boundary condition is recovered. No initial boundary compatibility is required. Since extension by zero can create a jump across \(\partial K\), the normal-trace estimate is stated away from \(t=0\); the density and interior flux estimates hold from \(t=0\). An auxiliary initialization removes this restriction for the trace estimate, as explained in Section~\ref{subsec:proof-strategy}. The result permits arbitrary interaction strength and does not require convexity of \(\Phi\) or \(W\).

The constant \(C_T\) may grow with \(T\). In the absence of interactinon, the dissipative structure of the linear equation prevents the approximation error from accumulating over time.

\begin{theorem}[Uniform-in-time convergence in the linear case]
\label{thm:linear-global}
Under the assumptions of Theorem~\ref{thm:strong-main}, suppose that \(W\equiv0\). There are \(C>0\) and \(\lambda_0>0\), independent of time, such that
\begin{align}
 &\sup_{t\ge0}\|\rho^\lambda(t)|_K-\rho(t)\|_{L^2(K)} \notag\\
 &\quad+\sup_{s\ge0}\left(\int_s^{s+1}
   \left[\|\rho^\lambda(t)|_K-\rho(t)\|_{H^1(K)}^2
   +\|J_K^\lambda(t)-J(t)\|_{L^2(K)}^2\right]\dd t\right)^{1/2}
 \le C\lambda^{1/4}
 \label{eq:linear-global-rate}
\end{align}
for \(0<\lambda\le\lambda_0\). Moreover,
\begin{equation}\label{eq:linear-global-L1}
  \sup_{t\ge0}
  \|\rho^\lambda(t)-\rho(t)\one_K\|_{L^1(\Rd)}
  \le C\lambda^{1/4}.
\end{equation}
\end{theorem}

Thus the \(L^2(K)\) density error has the same bound at every time, and the \(H^1(K)\) and flux errors have the same bound on every unit time interval. The \(L^1(\Rd)\) estimate also controls the exterior mass and gives the same rate in total variation and in the bounded--Lipschitz distance \(\dBL\), defined in \eqref{eq:BL-def}. Unit time intervals are natural here: for fixed \(\lambda>0\), the penalized and no-flux equilibrium densities differ even on \(K\), so their difference persists at large times.

Similar uniform-in-time estimates may also be available for interacting equations under suitable stability assumptions, for instance when the interaction is sufficiently weak. Such extensions are beyond the scope of the present work. The following example, however, shows that interaction can destroy uniform-in-time convergence when the no-flux dynamics starts at an unstable equilibrium. It uses a smooth cosine interaction and the same initial density for both equations.

\begin{theorem}[Failure of uniform-in-time convergence]
\label{thm:counterexample}
Let \(d=1\), \(K=[-\pi,\pi]\), \(\Phi=0\), and \(W(x,y)=-a\cos(x-y)\) with \(a>2\). Let \(\rho^\lambda\) be the penalized solution with initial density \(\bar\rho\one_K\), where \(\bar\rho=1/(2\pi)\), and let \(\rho\) be the corresponding no-flux solution. Then \(\rho(t)\equiv\bar\rho\). There exist \(\eta_a,C_a,\lambda_0>0\) such that, for every \(0<\lambda<\lambda_0\), one can find \(s_\lambda\) satisfying
\begin{equation}\label{eq:s-lambda-bound}
 2\le s_\lambda\le C_a(1+|\log\lambda|)
\end{equation}
and
\begin{equation}\label{eq:window-weak-main}
 \left|\int_{s_\lambda}^{s_\lambda+1}\!\int_K
 \sin^2x\,\bigl(\rho^\lambda(t,x)-\bar\rho\bigr)\dd x\dd t\right|
 \ge\eta_a.
\end{equation}
Consequently, convergence fails uniformly in time both narrowly on \(\R\) and weakly in \(L^2(K)\), and hence also strongly in \(L^2(K)\). Nevertheless, by Theorem~\ref{thm:strong-main}, \(\rho^\lambda|_K\to\bar\rho\) in \(L^\infty(0,T;L^2(K))\cap L^2(0,T;H^1(K))\) for every finite \(T>0\). In particular, \(s_\lambda\to\infty\) as \(\lambda\downarrow0\).
\end{theorem}

The discrepancy is detected by integration against a single fixed test function. It therefore reflects a loss of weak approximation, rather than only a loss of derivative control. Together, the three theorems distinguish approximation on fixed time intervals from approximation over the entire evolution. The latter holds without interaction, but requires further assumptions in the interacting case.

\subsection{Ideas of the proofs}\label{subsec:proof-strategy}

\paragraph{The boundary layer and the finite-time estimate.}
Inside \(K\), the penalization drift vanishes. The restriction of \(\rho^\lambda\) therefore satisfies the same type of equation as \(\rho\), with two errors: its normal flux need not vanish, and its interaction field includes the mass outside \(K\). Proposition~\ref{prop:error-estimate} reduces the comparison to these errors and the initial discrepancy. The exterior mass is \(O(\sqrt\lambda)\) on bounded time intervals. The main difficulty is to control the normal flux, which involves derivatives at the boundary and is not bounded by the exterior mass alone.

At distance \(s\) from \(K\), the confinement weight is \(w_\lambda=e^{-s^2/(2\lambda)}\). The exterior layer consequently has thickness \(O(\sqrt\lambda)\), while the restoring drift there has size \(O(\lambda^{-1/2})\). Passing to the relative density \(\rho^\lambda/w_\lambda\) absorbs this large drift into the weight. Weighted smoothing and tangential regularity then allow us to express the normal flux as an integral across the exterior layer; see Lemma~\ref{lem:transmission}. Cauchy--Schwarz in the normal variable extracts the square root of the layer's \(O(\sqrt\lambda)\) weight, giving the factor \(\lambda^{1/4}\) in Proposition~\ref{prop:residual}. This boundary estimate is the source of the rate in Theorem~\ref{thm:strong-main}.

To handle the jump in the zero-extended initial density, we first extend \(\rho_0\) across the boundary as a bounded \(H^1\) function, multiply it by \(w_\lambda\), and normalize its mass. The resulting auxiliary solution has finite initial weighted gradient energy and satisfies
\[
  \|\bar J_K^\lambda\cdot\vec n\|_{L^2(0,T;L^2(\partial K))}
  \le C_T\lambda^{1/4},
\]
where \(\bar J_K^\lambda\) is its interior flux. The auxiliary datum differs from the original one by \(O(\lambda^{1/4})\) in the weighted \(L^2\) norm. Stability of the penalized equation transfers the comparison to the common initial density; see Section~\ref{subsec:lifted} and Lemma~\ref{lem:pen-stability}.

\paragraph{Dissipation in the linear problem.}
When \(W=0\), the invariant densities are proportional to \(e^{-\Phi}\) on \(K\) and to \(e^{-\Phi-U_\lambda}\) on \(\Rd\). Lemma~\ref{lem:uniform-gibbs-gap} establishes a spectral gap uniform in \(\lambda\) for the latter, without assuming convexity of \(\Phi\). After subtracting equilibrium, the resulting dissipation controls the derivatives needed in the boundary-layer argument over the whole time interval \((0,\infty)\). This gives the global residual estimate in Lemma~\ref{lem:global-linear-residual}. A dissipative comparison on \(K\), proved in Lemma~\ref{lem:linear-dissipative}, then prevents the small boundary forcing from accumulating and yields Theorem~\ref{thm:linear-global}.

\paragraph{Amplification near an unstable equilibrium.}
For the cosine interaction in Theorem~\ref{thm:counterexample}, the no-flux solution remains exactly at the homogeneous equilibrium. Soft confinement creates an exterior layer, which generates a cosine moment of order \(\sqrt\lambda\). The corresponding linearized mode has growth rate \(a/2-1>0\), suggesting that this small perturbation becomes macroscopic within a time of order \(|\log\lambda|\). The proof turns this mechanism into an exact moment identity. Lemma~\ref{lem:boundary-forcing} identifies the initial forcing from the layer, and Proposition~\ref{prop:log-loss} shows that continued weak closeness on all unit intervals up to that time would force the cosine moment to exceed its bound for a probability density. This gives \eqref{eq:window-weak-main}.

\subsection{Background and related work}\label{subsec:background}

The replacement of a reflecting boundary by a confining drift has both PDE and probabilistic precedents. Alexander and Kim \cite{AlexanderKim2016} prove locally uniform approximation with rates for local quasilinear parabolic equations with oblique boundary data, using a cubic confining potential and a suitable extension of the operator. Alasio, Bruna, and Carrillo \cite{AlasioBrunaCarrillo2020} approximate no-flux solutions of nonlinear nonlocal Fokker--Planck equations by whole-space equations with increasingly strong confinement. Their weighted \(L^2\) and free-energy methods establish the limiting boundary condition and, in the linear-diffusion setting with interaction, strong space--time \(L^2\) convergence. For the quadratic Moreau--Yosida penalty considered here, our companion paper \cite{SeoKimChoi2026} proves qualitative convergence of the densities and interior fluxes and identifies the limiting boundary reaction. The present work obtains rates in density, gradient, and flux norms, and determines how the comparison can change on long time intervals.

At the probabilistic level, reflected diffusions and their construction were studied by Tanaka \cite{Tanaka1979}, Menaldi \cite{Menaldi1983}, and Lions and Sznitman \cite{LionsSznitman1984}. S{\l}omi\'nski established approximation results for reflected SDEs \cite{Slominski1994} and quantitative pathwise estimates for projection penalization with Lipschitz coefficients \cite{Slominski2013}. Bruggeman and Sarantsev \cite{BruggemanSarantsev2017} proved weak convergence for a class of penalty approximations on the half-line. These results concern paths or probability laws; the gradient and normal-flux estimates here require control of the density near the boundary. Nonlinear reflected diffusions and mean-field limits go back to Sznitman \cite{Sznitman1984}. Further developments include interaction through reflecting boundaries \cite{CoghiEtAl2022}, well-posedness and particle approximation on smooth domains \cite{HindsSharmaTretyakov2024}, and reflecting Fokker--Planck equations with nonlinear diffusion \cite{CiotirFayad2022}.

In constrained sampling, Brosse, Durmus, Moulines, and Pereyra \cite{BrosseDurmusMoulinesPereyra2017} use Moreau--Yosida penalization and estimate its bias for compactly supported log-concave targets. Projected Langevin sampling is studied in \cite{BubeckEldanLehec2018}, while penalized overdamped and underdamped methods with possibly nonconvex potentials are analyzed in \cite{GurbuzbalabanHuZhu2024}. Our estimates compare the corresponding continuous-time PDE evolutions from a common initial density, including their spatial gradients and fluxes.

Long-time behavior of interacting equations depends on the free-energy landscape. Bifurcation and phase transitions for periodic McKean--Vlasov equations, including the cosine interaction, are studied in \cite{CarrilloGvalaniPavliotisSchlichting2020}. Convergence to equilibrium and long-time particle approximation under additional structural assumptions have been obtained through logarithmic Sobolev inequalities \cite{Malrieu2001,Chizat2022} and, for kinetic equations, coupling estimates \cite{BolleyGuillinMalrieu2010}. Our counterexample shows how an unstable equilibrium can amplify the error created by soft confinement itself, despite convergence on every fixed time interval.

\paragraph{Organization.}
Section~\ref{sec:setting} gives the weak formulations and preliminary weighted estimates. Section~\ref{sec:strong} proves the boundary residual estimate and completes the finite-time comparison in Theorem~\ref{thm:strong-main}. Section~\ref{sec:linear-global} develops the uniform spectral-gap and dissipative estimates for Theorem~\ref{thm:linear-global}. Section~\ref{sec:counterexample} constructs the unstable example and proves Theorem~\ref{thm:counterexample}. The appendices supply the smoothing, normal-trace, geometric, and tangential regularity arguments.

\section{Setting and preliminaries}\label{sec:setting}
% ==================================================

The standing hypotheses stated in Section~\ref{sec:intro} remain in force. The preliminary results in this section do not require the \(C^3\) boundary regularity or the additional \(H^1\) initial regularity used in Theorems~\ref{thm:strong-main} and~\ref{thm:linear-global}. Unless stated otherwise, all $\lambda$-uniform estimates below are for $0<\lambda\le1$; fixed-parameter well-posedness is considered for every $\lambda>0$.
%The preliminary results in this section do not require the additional
%\(C^3\) boundary regularity or \(H^1\) initial regularity used in the first
%two theorems.

\subsection{Equations and notation}\label{subsec:notation}

Fix a time horizon \(T>0\). Write
\[
  \Gamma=\partial K,\qquad
  U_\lambda(x)=\frac{\dist(x,K)^2}{2\lambda},\qquad
  w_\lambda=e^{-U_\lambda},\qquad
  \dd\widehat\mu_\lambda=w_\lambda\,\dd x .
\]
The measure \(\widehat\mu_\lambda\) is not normalized and \(w_\lambda=1\) on \(K\). We write $\vec n(y)$ for the outward unit normal at $y\in\Gamma$ and extend this notation to $K^c$ by $\vec n(x):=\nabla\dist(x,K)$. For $\mu\in\mcP(\Rd)$, recall the drift field
\begin{equation*}
  b[\mu](x):=\nabla\Phi(x)+\int_{\Rd}\nabla_xW(x,y)\dd\mu(y).
\end{equation*}
The \emph{penalized interacting dynamics} is the McKean--Vlasov Fokker--Planck equation
\begin{equation}\label{eq:pen-pde}
  \partial_t\rho^\lambda=\divg J^\lambda,
  \qquad
  J^\lambda:=\nabla\rho^\lambda+\rho^\lambda\nabla U_\lambda
  +\rho^\lambda\,b[\rho^\lambda(t)],
  \qquad\text{in }(0,T)\times\Rd,
\end{equation}
with initial datum $\rho^\lambda(0,\cdot)=\rho_0$. The \emph{no-flux interacting dynamics} is
\begin{equation}\label{eq:ref-pde}
\begin{cases}
  \partial_t\rho=\divg J, & \text{in }(0,T)\times K,\\[0.2em]
  J\cdot\vec n=0, & \text{on }(0,T)\times\Gamma,\\[0.2em]
  \rho(0,\cdot)=\rho_0, & \text{in }K,
\end{cases}
\qquad
  J:=\nabla\rho+\rho\,b[\rho(t)].
\end{equation}
We distinguish the \emph{full flux} $J^\lambda$ of \eqref{eq:pen-pde} from the \emph{interior (unpenalized) flux} $J_K^\lambda$ in \eqref{eq:JK-def}. They coincide on $K$, where $\nabla U_\lambda=0$, but not outside. The no-flux condition in \eqref{eq:ref-pde} concerns the limit of $J_K^\lambda\cdot\vec n$.

Since $K$ is closed and convex, $\proj_K$ is single-valued and $1$-Lipschitz on $\Rd$, and $\dist^2\in C^{1,1}(\Rd)$. The $C^2$ regularity of $\partial K$ further implies that $\dist$ is $C^2$ on all of $K^c$,  with
\begin{equation*}
  0\le\Delta\dist(x)\le C_K,\qquad
  |D^2\dist(x)|\le C_K\quad\text{for every }x\in K^c,
\end{equation*}
where $C_K$ depends only on the curvature bounds of $\partial K$. We will also use the Gaussian layer estimates of Lemma~\ref{lem:layer}.

For the drift estimates, set
\[
  B:=\max\Bigl\{
    \|\nabla\Phi\|_\infty+\|\nabla_xW\|_\infty,
    \|\nabla^2\Phi\|_\infty+\|\nabla_x^2W\|_\infty
  \Bigr\},\qquad
  L_W:=\|\nabla_xW\|_\infty .
\]
Then, for all $\mu,\nu\in\mcP(\Rd)$ with densities $\rho_\mu,\rho_\nu\in L^1$,
\begin{equation}\label{eq:drift-bounds}
  \|b[\mu]\|_\infty\le B,\qquad
  \|\nabla b[\mu]\|_\infty\le B,\qquad
  \|b[\mu]-b[\nu]\|_\infty\le L_W\,\|\rho_\mu-\rho_\nu\|_{L^1(\Rd)} .
\end{equation}
The drift class
\begin{equation}\label{eq:drift-class}
  \mathcal B:=\bigl\{\,b\in C^1(\Rd;\Rd):\ \|b\|_\infty\le B,\
  \|\nabla b\|_\infty\le B\,\bigr\}
\end{equation}
contains $b[\mu]$ for every probability measure $\mu$, uniformly in $\lambda$ and $t$.

We write $M:=\|\rho_0\|_{L^\infty(K)}$ and, when $\rho_0\in H^1(K)$, $N_0:=1+M+\|\rho_0\|_{H^1(K)}$.

In the finite-time estimates, $C_T$ denotes a positive constant depending on $T$, $K$, $d$, $B$, $\|W\|_{C^2}$, $\|\rho_0\|_{L^\infty(K)}$, and $\|\rho_0\|_{L^2(K)}$; its value may change from line to line. Dependence on the $H^1$ norm of the initial datum is displayed through $N_0$ in the proofs and may be absorbed into the constants in Theorems~\ref{thm:strong-main} and \ref{thm:linear-global}.

\subsection{Weak formulations}\label{subsec:weak-formulations}

We use a weighted formulation for the penalized equation and the standard variational formulation for the no-flux equation.

\begin{definition}\label{def:pen-sol}
A \emph{penalized solution} on $[0,T]$ is a nonnegative \(\rho^\lambda\in C([0,T];L^1(\Rd))\) with \(\rho^\lambda(t)\in\mcP(\Rd)\). Its relative density \(\hh^\lambda:=\rho^\lambda/\wl\) satisfies
\[
 \hh^\lambda\in
 L^\infty(0,T;L^2(\muh))\cap L^2(0,T;H^1(\muh)),
\]
and, for every \(\varphi\in C_c^\infty([0,T)\times\Rd)\),
\begin{equation}\label{eq:pen-weak}
  \int_{\Rd}\rho_0\varphi(0)\dd x
  +\int_0^T\!\!\int_{\Rd}\rho^\lambda\,\partial_t\varphi\dd x\dd t
  =\int_0^T\!\!\int_{\Rd}
  \bigl(\wl\nabla\hh^\lambda+\rho^\lambda\,b[\rho^\lambda(t)]\bigr)
  \cdot\nabla\varphi\dd x\dd t,
\end{equation}
where $H^1(\muh):=\{v:\int(v^2+|\nabla v|^2)\dd\muh<\infty\}$.
\end{definition}

Since \(\wl\nabla\hh^\lambda=\nabla\rho^\lambda+\rho^\lambda\nabla U_\lambda\), \eqref{eq:pen-weak} is the weak form of \eqref{eq:pen-pde}, with \(J^\lambda=\wl\nabla\hh^\lambda+\rho^\lambda b[\rho^\lambda(t)]\). This identity is the reason for working with \(\hh^\lambda\). When \(b\equiv0\), every profile \(c\wl\), with \(c\) constant, has zero flux, \(\nabla(c\wl)+c\wl\nabla U_\lambda=0\), although at distance of order \(\sqrt\lambda\) from \(K\) each of the two terms is of order \(\lambda^{-1/2}c\wl\), as \(|\nabla U_\lambda|=\dist(\cdot,K)/\lambda\). Writing \(\rho^\lambda=\wl\hh^\lambda\) factors out this profile, and the identity replaces the sum \(\nabla\rho^\lambda+\rho^\lambda\nabla U_\lambda\) by the single term \(\wl\nabla\hh^\lambda\); the cancellation therefore need not be recovered from separate estimates of the two terms. Accordingly, the whole-space energy estimates for the penalized flow in this section and in Section~\ref{sec:strong} are formulated for \(\hh^\lambda\) in \(L^2(\muh)\) and \(H^1(\muh)\). Because \(\wl=1\) and \(\nabla\wl=0\) on \(K\), they also control \(\rho^\lambda|_K\) in \(L^2(K)\) and \(H^1(K)\).

For an auxiliary flow, the same definition is used with its stated initial probability density in place of $\rho_0$.

\begin{definition}
\label{def:weak-sol}
A \emph{weak solution} of \eqref{eq:ref-pde} is a nonnegative
\[
  \rho\in L^\infty(0,T;L^2(K))\cap L^2(0,T;H^1(K)),\qquad
  \partial_t\rho\in L^2(0,T;H^1(K)'),
\]
with $\rho(0)=\rho_0$, such that for a.e.\ $t$ and all $\varphi\in H^1(K)$,
\begin{equation}\label{eq:weak-form}
  \langle\partial_t\rho(t),\varphi\rangle
  +\int_K\bigl(\nabla\rho+\rho\,b[\rho(t)]\bigr)\cdot\nabla\varphi\dd x=0 .
\end{equation}
No boundary integral appears in \eqref{eq:weak-form}; this encodes $J\cdot\vec n=0$ in the normal-trace sense. Weak solutions belong to $C([0,T];L^2(K))$ and conserve mass.
\end{definition}

For finite Borel measures \(\mu,\nu\) on \(\Rd\), we use the bounded--Lipschitz distance
\begin{equation}\label{eq:BL-def}
  \dBL(\mu,\nu):=
  \sup_{\|\varphi\|_\infty+\operatorname{Lip}(\varphi)\le1}
  \left|\int_{\Rd}\varphi\,\dd(\mu-\nu)\right|.
\end{equation}
On probability measures, this distance metrizes narrow convergence: \(\dBL(\mu_n,\mu)\to0\) if and only if \(\mu_n\to\mu\) narrowly.

\subsection{Results from the companion paper}\label{subsec:companion-inputs}

We recall the well-posedness and qualitative limit from \cite{SeoKimChoi2026}, together with the weighted bounds used below. The qualitative convergence identifies the limit; the strong rates are proved in Section~\ref{sec:strong}. Weak solutions in the sense of Definition~\ref{def:weak-sol} are called variational solutions in \cite[Definition~2.8]{SeoKimChoi2026}.
% Companion locators checked against arXiv:2610.25025v1 (20 October 2026):
% Definition 2.8, Theorem 2.9, Propositions 3.1--3.2, Theorem 4.1,
% Corollary 4.2, Lemma 4.6, and Lemma C.1. Recheck these locators
% when citing a later version.

The following qualitative limit follows from \cite[Theorem~2.9 and Propositions~3.1--3.2]{SeoKimChoi2026}.

\begin{theorem}
\label{thm:weak-main}
For each $\lambda>0$ the penalized equation has a unique solution, and the no-flux equation has a unique weak solution $\rho$. For every $T<\infty$, as $\lambda\downarrow0$,
\begin{gather*}
  \rho^\lambda|_K\longrightarrow\rho
  \quad\text{strongly in }L^2((0,T)\times K),\\
  \sup_{0\le t\le T}
  \dBL(\rho^\lambda(t),\rho(t)\one_K)\longrightarrow0,\\
  J^\lambda|_K\rightharpoonup J
  \quad\text{weakly in }L^2((0,T)\times K),
\end{gather*}
where \(\rho\) is the unique weak solution of \eqref{eq:ref-pde}. Moreover,
\begin{equation}\label{eq:prior-exterior}
  \sup_{0\le t\le T}\int_{K^c}\rho^\lambda(t)\,\dd x
  \le C_T\sqrt\lambda .
\end{equation}
\end{theorem}

Moreover, the fixed-parameter well-posedness result \cite[Proposition~3.1]{SeoKimChoi2026} includes the weighted energy identity
\begin{equation}\label{eq:energy-identity}
  \frac12\frac{\dd}{\dd t}\int_{\Rd}|\hh^\lambda(t)|^2\,\dd\muh
  +\int_{\Rd}|\nabla\hh^\lambda(t)|^2\,\dd\muh
  =-\int_{\Rd}\hh^\lambda(t)\,b[\rho^\lambda(t)]\cdot
    \nabla\hh^\lambda(t)\,\dd\muh
\end{equation}
for a.e.\ $t>0$. All energy identities below are justified by Galerkin approximation or Steklov averaging; we use the Lions--Magenes chain rule in the Gelfand triple $H^1(K)\subset L^2(K)\subset H^1(K)'$.

The auxiliary flows in Sections~\ref{subsec:lifted} and \ref{subsec:linear-lifted} use initial probability densities $h_0w_\lambda$ with $h_0\in L^2(\muh)$, rather than data supported in $K$. These data are explicitly covered by \cite[Proposition~3.1]{SeoKimChoi2026}, which gives well-posedness and the energy identity \eqref{eq:energy-identity} for initial probability densities with relative density in \(L^2(\muh)\).

We use the following drift-uniform domination estimate from \cite[Theorem~4.1, Corollary~4.2, and Lemma~C.1]{SeoKimChoi2026} in the form needed for both initializations.
\begin{proposition}
\label{thm:barrier}
There are constants $\beta,\gamma>0$, depending only on $K,d$, and the drift class \eqref{eq:drift-class}, with the following property. Let \(\rho^\lambda\) be a penalized solution whose initial datum satisfies $0\le\rho_0^\lambda\le M_0e^{-\dist^2/(2\lambda)}$. Then
\begin{equation}\label{eq:domination}
  0\le\rho^\lambda(t,x)
  \le C M_0e^{\gamma t}
  \exp\!\left(-\frac{\dist(x)^2}{2\lambda}+\beta\dist(x)\right).
\end{equation}
Moreover, for every $T<\infty$ and $q\ge0$,
\begin{equation}\label{eq:exterior-moments}
  \int_{K^c}\dist(x)^q\rho^\lambda(t,x)\dd x
  \le C_{T,q}\lambda^{(q+1)/2},\qquad
  \int_{\{\dist\ge\delta\}}\rho^\lambda(t,x)\dd x
  \le C_Te^{-\delta^2/(8\lambda)}
\end{equation}
for $t\in[0,T]$ and $0<\delta\le1$, and $\|\rho^\lambda\|_{L^\infty([0,T]\times K)}\le C_TM_0$. These conclusions apply to both the original and auxiliary initial densities below.
\end{proposition}

In \cite[Theorem~4.1]{SeoKimChoi2026}, the function \(M_0e^{\gamma t}e^{\beta\zeta}\wl\) is shown to be a supersolution of the linear penalized equation for every measurable family of drifts in the class used there, which contains \(b[\mu]\) for all \(\mu\in\mcP(\Rd)\); the tilt \(\zeta\ge0\) is bounded on \(K\), and \(\zeta\le\dist(\cdot,K)+C_K\) on \(K^c\). That theorem explicitly includes initial data \(0\le\rho_0^\lambda\le M_0\wl\), and therefore gives \eqref{eq:domination} for both initializations. As in \cite[Corollary~4.2]{SeoKimChoi2026}, \eqref{eq:exterior-moments} and the interior bound follow from \eqref{eq:domination} and Lemma~\ref{lem:layer}, for every \(0<\lambda\le1\); the case \(q=0\) gives \eqref{eq:prior-exterior} in the same range.

\subsection{Basic weighted estimates}\label{subsec:basic-estimates}

For the weighted estimates below $\rho^\lambda$ is the penalized solution, $\hh:=\hh^\lambda=\rho^\lambda/\wl$ (so $\hh=\rho^\lambda$ on $K$), and $b^\lambda(t,x):=b[\rho^\lambda(t)](x)\in\mathcal B$. In terms of $\hh$ the equation reads, weakly,
\begin{equation}\label{eq:ratio-eq}
  \wl\,\partial_t\hh=\divg\bigl(\wl\nabla\hh+\wl\,\hh\,b^\lambda\bigr),
\end{equation}
because $\nabla\rho^\lambda+\rho^\lambda\nabla U_\lambda=\wl\nabla\hh$. We write
\[
\begin{aligned}
  \|v\|_\lambda^2
  &:=\int_{\Rd}v^2\dd\muh
   =\int_{\Rd}v^2w_\lambda\,\dd x,\\
  E(t)
  &:=\int_{\Rd}|\nabla\hh(t)|^2\dd\muh
   =\int_{\Rd}|\nabla\hh(t)|^2w_\lambda\,\dd x.
\end{aligned}
\]
Recall $\muh(\Rd)\le|K|+C\sqrt\lambda\le C_K$ and $\|\hh_0\|_\lambda=\|\rho_0\|_{L^2(K)}$ for the hard initialization. The energy, drift-regularity, smoothing, and residual estimates also apply to initial probability densities $\hh_0\wl$ with $\hh_0\in L^2(\muh)$. For these weighted-data estimates, $C_T$ may depend on a $\lambda$-independent upper bound for $\|\hh_0\|_\lambda$; dependence on the initial energy is displayed explicitly.

\begin{lemma}\label{lem:basic-energy}
For all $t\in[0,T]$,
\[
  \|\hh(t)\|_\lambda^2+\int_0^tE(s)\dd s
  \ \le\ e^{B^2t}\,\|\hh_0\|_\lambda^2 .
\]
\end{lemma}

\begin{proof}
The identity \eqref{eq:energy-identity} and Cauchy--Schwarz give $\frac12\frac{\dd}{\dd t}\|\hh\|_\lambda^2+E \le B\|\hh\|_\lambda E^{1/2} \le\frac12E+\frac{B^2}2\|\hh\|_\lambda^2$. Gr\"onwall's inequality yields the assertion.
\end{proof}

\begin{lemma}
\label{lem:drift-time-reg}
Set
\[
  \dot b^\lambda(t,x):=-\int_{\Rd}\nabla_y\nabla_xW(x,y)\cdot
  J^\lambda(t,y)\dd y .
\]
Then $t\mapsto b^\lambda(t,x)$ is absolutely continuous for every $x$, with
\begin{equation}\label{eq:drift-ac}
  b^\lambda(t,x)-b^\lambda(s,x)
  =\int_s^t \dot b^\lambda(r,x)\dd r ;
\end{equation}
moreover,
\begin{equation}\label{eq:J-L1-bound}
  \sup_{x\in\Rd}|\dot b^\lambda(t,x)|
  \le \|W\|_{C^2}\int_{\Rd}|J^\lambda(t,y)|\dd y
  \le C\bigl(1+E(t)^{1/2}\bigr).
\end{equation}
In particular, the left-hand side of \eqref{eq:J-L1-bound}, regarded as a function of $t$, is bounded in $L^2(0,T)$ uniformly in $\lambda$.
\end{lemma}

\begin{proof}
The full-flux bound of \cite[Lemma~4.6]{SeoKimChoi2026} also follows directly from $J^\lambda=w_\lambda\nabla\hh+\rho^\lambda b^\lambda$:
\[
  \int_{\Rd}|J^\lambda(t)|\dd x
  \le \muh(\Rd)^{1/2}E(t)^{1/2}+B.
\]
For fixed $x$, test the continuity equation with smooth compactly supported approximations of $y\mapsto\nabla_xW(x,y)$. The bounded first and mixed second derivatives of $W$, together with the preceding $L^1$ bound, permit passage to the limit in the spatial cutoff and the smooth approximation. This gives \eqref{eq:drift-ac} in the sense of distributions in time. Its right-hand side belongs to $L^2(0,T)$ by Lemma~\ref{lem:basic-energy}, which proves absolute continuity and \eqref{eq:J-L1-bound}.
\end{proof}

\section{Proof of Theorem~\ref{thm:strong-main}}\label{sec:strong}

We work under the hypotheses of Theorem~\ref{thm:strong-main}. The comparison on \(K\) is carried out in Proposition~\ref{prop:error-estimate}: for a penalized solution whose \(L^2(K)\)-norm is bounded and whose exterior mass is \(O(\sqrt\lambda)\), uniformly on \([0,T]\), its estimate \eqref{eq:error-estimate} bounds the squared density and gradient errors on \(K\) by \(C_T\) times the sum of the squared initial mismatch, a term \(\lambda\) coming from the exterior mass in the interaction drift, and the squared \(L^2((0,T)\times\Gamma)\)-norm of the residual normal flux. The main task is therefore to show that this residual is \(O(\lambda^{1/4})\). Since the conormal flux is continuous across \(\Gamma\) (Lemma~\ref{lem:transmission}), integrating the collar equation over the exterior side of the collar expresses the residual as an integral in the normal variable \(s\) against \(e^{-s^2/(2\lambda)}\), up to a bounded Jacobian; this factor has integral \(O(\sqrt\lambda)\) over \(s>0\) and is independent of the tangential variables. Weighted smoothing (Section~\ref{sec:energy}) and tangential regularity (Section~\ref{subsec:collar}) justify this representation and control the integrand in weighted \(L^2\), uniformly in \(\lambda\); Cauchy--Schwarz in \(s\) then yields the factor \(\lambda^{1/4}\) (Proposition~\ref{prop:residual}).

For the hard datum, which jumps across \(\Gamma\) unless \(\rho_0|_\Gamma\equiv0\), this residual bound is available only on \((\tau,T)\), with a factor \(1+\tau^{-1/2}\). We therefore apply Proposition~\ref{prop:error-estimate} to an auxiliary solution with finite initial weighted gradient energy, whose residual bound holds down to \(t=0\), and return to the hard solution by \(L^2(\muh)\)-stability of the penalized flow; the initial mismatch is also \(O(\lambda^{1/4})\) (Section~\ref{subsec:lifted}). The interior flux error then follows from the resulting density estimates and the exterior-mass bound, and \eqref{eq:bdry-rate} is Proposition~\ref{prop:residual} applied directly to the hard solution.

\subsection{Weighted smoothing and interior bounds}\label{sec:energy}

We use the relative density, weighted norm, and energy introduced in Section~\ref{subsec:basic-estimates}.

\begin{lemma}\label{lem:smoothing}
There is a constant $C_T$, independent of $\lambda$, such that:
\begin{enumerate}[label=(\roman*),leftmargin=2em]
\item ($L^2$ data) for every $\tau\in(0,T]$,
\[
  \sup_{\tau\le t\le T}E(t)
  +\int_\tau^T\|\partial_t\hh(s)\|_\lambda^2\dd s
  \ \le\ C_T\,(1+\tau^{-1});
\]
\item ($H^1$ data) if moreover $E(0)<\infty$, then
\[
  \sup_{0\le t\le T}E(t)
  +\int_0^T\|\partial_t\hh(s)\|_\lambda^2\dd s
  \ \le\ C_T\,\bigl(1+E(0)\bigr).
\]
\end{enumerate}
\end{lemma}

\begin{proof}
Set
\[
  Y(t):=\int_{\Rd}\hh\,b^\lambda\cdot\nabla\hh\,\dd\muh.
\]
Testing the equation by \(\partial_t\hh\) and applying the product rule to \(Y\) formally gives
\[
\begin{aligned}
  \|\partial_t\hh\|_\lambda^2+\frac12E'+Y'
  &=\int_{\Rd}\partial_t\hh\,b^\lambda\cdot\nabla\hh\,\dd\muh
    +\int_{\Rd}\hh\,\dot b^\lambda\cdot\nabla\hh\,\dd\muh\\
  &\le\frac12\|\partial_t\hh\|_\lambda^2+\frac{B^2}{2}E
    +\|\dot b^\lambda\|_\infty\|\hh\|_\lambda E^{1/2}.
\end{aligned}
\]
Thus no estimate of \(\nabla\partial_t\hh\) is needed. Moreover,
\[
  |Y(t)|\le B\|\hh(t)\|_\lambda E(t)^{1/2}
  \le\frac14E(t)+B^2\|\hh(t)\|_\lambda^2.
\]
The basic energy estimate and the \(L^2_tL^\infty_x\) bound for \(\dot b^\lambda\) from Lemma~\ref{lem:drift-time-reg} control the time integral of the last term. Integrating and using the bound on \(Y\) at the endpoints gives, for \(H^1(\muh)\) initial data,
\[
  \int_0^T\|\partial_t\hh(t)\|_\lambda^2\,\dd t
  +\sup_{0\le t\le T}E(t)\le C_T(1+E(0)).
\]
For \(L^2(\muh)\) initial data, the same calculation with \(t\partial_t\hh\) gives
\[
  \int_0^Tt\|\partial_t\hh(t)\|_\lambda^2\,\dd t
  +\sup_{0<t\le T}tE(t)\le C_T.
\]
These inequalities imply (ii) and (i), respectively. The Galerkin approximation, the product rule for the time-dependent drift term, and the identification of the limit are given in Appendix~\ref{app:smoothing}.
\end{proof}

\begin{lemma}
\label{lem:interior}
There is a constant $C_T$, independent of $\lambda$, such that
\begin{align}
  &\sup_{t\le T}\|\rho^\lambda(t)\|_{L^2(K)}
  +\|\nabla\rho^\lambda\|_{L^2((0,T)\times K)}
  +\|J^\lambda\|_{L^2((0,T)\times K)}\le C_T,
  \label{eq:interior-bounds}\\
  &\|\partial_t\rho^\lambda\|_{L^2(0,T;H^1_0(K)')}\le C_T .
  \label{eq:dual-bound}
\end{align}
\end{lemma}

\begin{proof}
On $K$, $\wl\equiv1$, so \eqref{eq:interior-bounds} is Lemma~\ref{lem:basic-energy} together with $|J^\lambda|\le|\nabla\hh|+B|\hh|$ on $K$ (the penalty drift vanishes on $K$). To prove \eqref{eq:dual-bound}, we test \eqref{eq:pen-weak} with $\varphi\in C_c^\infty(K)$. Then $\langle\partial_t\rho^\lambda,\varphi\rangle =-\int_KJ^\lambda\cdot\nabla\varphi$, and \eqref{eq:interior-bounds} gives \eqref{eq:dual-bound}.
\end{proof}

\subsection{Collar coordinates and tangential regularity}\label{subsec:collar}

Let $\Gamma:=\partial K$ and, for $\delta\in(0,\delta_0]$ fixed (depending only on $K$), let
\[
  \Psi:\Gamma\times(-\delta,\delta)\to
  \{|d_s|<\delta\},\qquad
  \Psi(y,s):=y+s\,\vec n(y),
\]
which is a $C^2$ diffeomorphism since $\partial K\in C^3$; here $d_s$ is the signed distance ($d_s<0$ in $K$, $d_s=\dist$ outside), and $d_s(\Psi(y,s))=s$. In these coordinates the weight is
\[
  \wl(\Psi(y,s))=e^{-s_+^2/2\lambda},
\]
a function of $s$ alone. Because the weight depends only on $s$, tangential differentiation introduces no $\lambda$-dependent factor. Let $m(y,s)$ denote the Jacobian of $\Psi$ ($m\in C^1$, $m(y,0)=1$, $c\le m\le C$ on the collar), so that $\dd x=m(y,s)\dd\Hm(y)\dd s$. Set
\[
  a_\lambda(y,s):=e^{-s_+^2/2\lambda}\,m(y,s),
\]
so that $\dd\muh=a_\lambda\dd\Hm(y)\dd s$ on the collar. Writing
\[
q(t,y,s):=\hh(t,\Psi(y,s)),
\] 
the weak equation \eqref{eq:ratio-eq} transforms into
\begin{equation}\label{eq:collar-eq}
  a_\lambda\,\partial_tq
  =\partial_s\Bigl(a_\lambda\bigl(\partial_sq+\beta_s\,q\bigr)\Bigr)
  +\divG\Bigl(a_\lambda\bigl(\mathsf B\,\nG q+\beta_\Gamma\,q\bigr)\Bigr)
  \qquad\text{on }\Gamma\times(-\delta,\delta),
\end{equation}
in the weak sense against test functions supported in the collar. To specify the coefficients, let $S_y:=D_\Gamma\vec n(y)$ and $A_s(y):=I+sS_y$ on $T_y\Gamma$, and let $P_\Gamma(y)$ be the orthogonal projection onto $T_y\Gamma$. The tangential differential of $\Psi$ is $A_s$, so $m=\det A_s$ and $\nabla\hh\circ\Psi=\vec n\,\partial_sq+A_s^{-1}\nG q$. Thus the coefficients in \eqref{eq:collar-eq} are
\[
 \begin{aligned}
 \mathsf B(y,s)&=A_s(y)^{-2},\\
 \beta_s(t,y,s)&=b^\lambda(t,\Psi(y,s))\cdot\vec n(y),\\
 \beta_\Gamma(t,y,s)&=A_s(y)^{-1}P_\Gamma(y)b^\lambda(t,\Psi(y,s)).
 \end{aligned}
\]
In particular, $\mathsf B$ is symmetric and uniformly elliptic on the collar. By the drift bounds \eqref{eq:drift-bounds} and the $C^2$ regularity of $\Psi$,
\begin{equation}\label{eq:collar-coeff-bounds}
  \|\beta_s\|_{L^\infty_t C^1_{y,s}}
  +\|\beta_\Gamma\|_{L^\infty_t C^1_{y,s}}
  +\|\mathsf B\|_{C^1}
  +\|\nG\log m\|_\infty+\|\partial_s\log m\|_\infty\ \le\ C_K(1+B),
\end{equation}
uniformly in $t$ and $\lambda$, and $\|\partial_t\beta_s\|_{L^2_tL^\infty}+\|\partial_t\beta_\Gamma\|_{L^2_tL^\infty} \le C$ by Lemma~\ref{lem:drift-time-reg}. Moreover, $\nG\log a_\lambda=\nG\log m$ is bounded independently of $\lambda$: the singular factor is never differentiated tangentially.

Let us define the \emph{conormal flux profile}
\begin{equation*}
  F(t,y,s):=a_\lambda(y,s)\,\bigl(\partial_sq+\beta_s\,q\bigr)(t,y,s),
\end{equation*}
so that \eqref{eq:collar-eq} reads $a_\lambda\partial_tq=\partial_sF+ \divG\bigl(a_\lambda(\mathsf B\nG q+\beta_\Gamma q)\bigr)$. On the interior side ($s<0$), where $\wl\equiv1$ and $\rho^\lambda=q\circ\Psi^{-1}$, the chain rule gives
\begin{equation*}
  \bigl(J^\lambda_K\cdot\vec n\bigr)\bigl(t,\Psi(y,s)\bigr)
  =\bigl(\partial_sq+\beta_sq\bigr)(t,y,s)
  =\frac{F(t,y,s)}{m(y,s)},
  \qquad s<0 ,
\end{equation*}
for a.e.\ $(y,s)$.

\paragraph{Tangential regularity.}\label{subsec:tangential}

For a.e.\ fixed $t$, we regard \eqref{eq:collar-eq} as an elliptic equation for $q(t)$, with $f:=\partial_tq(t)\in L^2(a_\lambda)$ treated as given data:
\[
  -\partial_s\bigl(a_\lambda(\partial_sq+\beta_sq)\bigr)
  -\divG\bigl(a_\lambda(\mathsf B\nG q+\beta_\Gamma q)\bigr)
  =-a_\lambda f.
\]

\begin{lemma}\label{lem:tangential}
There is $C_T>0$ such that, for almost every $t\in(0,T)$,
\[
  \bigl\|\nG^2q(t)\bigr\|_{L^2(a_\lambda)}^2
  +\bigl\|\partial_s\nG q(t)\bigr\|_{L^2(a_\lambda)}^2
  \ \le\
  C_T\Bigl(\|\partial_t\hh(t)\|_\lambda^2+E(t)+\|\hh(t)\|_\lambda^2\Bigr),
\]
where $\|\cdot\|_{L^2(a_\lambda)}$ is the norm on $\Gamma\times(-\delta',\delta')$, $\delta'=\delta/2$, with measure $a_\lambda\dd\Hm\dd s$, and $\nG^2$, $\partial_s\nG$ denote second derivatives in a fixed finite atlas.
\end{lemma}

The proof by tangential difference quotients is given in Appendix~\ref{app:dq}. Tangential translations leave the singular factor $e^{-s_+^2/2\lambda}$ unchanged, and $|D^\eta a_\lambda|\le Ca_\lambda$ uniformly in $\lambda$.

\subsection{The normal trace and the boundary residual}
\label{subsec:transmission}

An energy bound alone does not give an $L^2(\Gamma)$ normal trace. The next lemma uses tangential regularity to make the conormal profile continuous across $s=0$ and identify its value with that trace.

\begin{lemma}\label{lem:transmission}
For almost every $t\in(0,T)$, after modification on a null set of $s$,
\[
  F(t,\cdot,\cdot)\in C\bigl([-\delta',\delta'];L^2(\Gamma)\bigr),
\]
with distributional $s$-derivative
\begin{equation}\label{eq:dsF}
  \partial_sF=a_\lambda\,\partial_tq
  -\divG\bigl(a_\lambda(\mathsf B\nG q+\beta_\Gamma q)\bigr)
  \ \in\ L^1\bigl((-\delta',\delta');L^2(\Gamma)\bigr).
\end{equation}
Moreover, the value $F(t,\cdot,0)\in L^2(\Gamma)$, taken from this continuous representative, is the normal trace of the interior flux:
\begin{equation}\label{eq:g-def}
  g^\lambda(t,\cdot):=F(t,\cdot,0)
  =\lim_{\eps\downarrow0}
  \bigl(J^\lambda_K\cdot\vec n\bigr)\bigl(t,\Psi(\cdot,-\eps)\bigr)
  \quad\text{in }L^2(\Gamma),
\end{equation}
and for every $\varphi\in H^1(K)$,
\begin{equation}\label{eq:trace-green}
  \int_K J^\lambda_K(t)\cdot\nabla\varphi\dd x
  +\int_K\bigl(\divg J^\lambda_K(t)\bigr)\varphi\dd x
  =\int_\Gamma g^\lambda(t,y)\,\varphi(y)\dd\Hm(y),
\end{equation}
where $\divg J^\lambda_K(t)=\partial_t\rho^\lambda(t)\in L^2(K)$.
\end{lemma}

\begin{proof}
Fix a time \(t\) for which the fixed-time weak formulation, \(\partial_t\hh(t)\in L^2(\muh)\), and Lemma~\ref{lem:tangential} hold; these conditions hold for almost every \(t\in(0,T)\). We suppress \(t\) from the notation and put \(I=(-\delta',\delta')\).

For every measurable \(v\),
\begin{equation}\label{eq:weighted-L1-slice}
  \int_I\|a_\lambda(\cdot,s)v(\cdot,s)\|_{L^2(\Gamma)}\,\dd s
  \le C\|v\|_{L^2(a_\lambda\,\dd\Hm\,\dd s)}.
\end{equation}
Indeed, \(\|a_\lambda v\|_{L^2(\Gamma)}^2 \le\|a_\lambda(\cdot,s)\|_{L^\infty(\Gamma)} \int_\Gamma a_\lambda v^2\,\dd\Hm\), and \(\int_I\|a_\lambda(\cdot,s)\|_{L^\infty(\Gamma)}\,\dd s\le C\).

Recall that \(F=a_\lambda(\partial_sq+\beta_sq)\). The energy estimate and \eqref{eq:weighted-L1-slice} imply \(F\in L^1(I;L^2(\Gamma))\). Since $\partial_t\hh(t)\in L^2(\muh)$, the change of variables to collar coordinates gives $\partial_tq\in L^2(a_\lambda\dd\Hm\dd s)$; hence \eqref{eq:weighted-L1-slice}, applied with $v=\partial_tq$, yields $a_\lambda\partial_tq\in L^1(I;L^2(\Gamma))$. In every boundary chart,
\begin{equation}\label{eq:tangential-div-bound}
  \left|\divG\bigl(a_\lambda(\mathsf B\nG q+\beta_\Gamma q)\bigr)\right|
  \le Ca_\lambda\bigl(|\nG^2q|+|\nG q|+|q|\bigr).
\end{equation}
No normal derivative of the Gaussian factor occurs here: \(\nG a_\lambda=a_\lambda\nG\log m\). Thus Lemma~\ref{lem:tangential} and \eqref{eq:weighted-L1-slice} also give $\divG\bigl(a_\lambda(\mathsf B\nG q+\beta_\Gamma q)\bigr) \in L^1(I;L^2(\Gamma))$. Testing the fixed-time collar equation against tensor products and using density, we obtain, for every $\eta\in C_c^\infty(I)$,
\[
  -\int_I F(s)\eta'(s)\dd s
  =\int_I\bigl(a_\lambda\partial_tq
    -\divG\bigl(a_\lambda(\mathsf B\nG q+\beta_\Gamma q)\bigr)\bigr)(s)
    \eta(s)\dd s
  \quad\text{in }L^2(\Gamma).
\]
Thus the $L^2(\Gamma)$-valued distributional derivative of $F$ belongs to $L^1(I;L^2(\Gamma))$, and consequently \(F\in W^{1,1}(I;L^2(\Gamma))\hookrightarrow C(\overline I;L^2(\Gamma))\). In particular the one-sided values at \(s=0\) coincide. Define \(g^\lambda(t)=F(t,\cdot,0)\).

On the interior side, \(F(y,s)=m(y,s) (J^\lambda_K\cdot\vec n_s)(\Psi(y,s))\) for \(s<0\), where \(\vec n_s(\Psi(y,s))=\vec n(y)\). Its \(L^2(\Gamma)\)-continuity therefore gives the parallel-surface limit in \eqref{eq:g-def}. The identification of this limit with the \(H(\divg;K)\) normal trace, including the Green identity \eqref{eq:trace-green}, is proved in Appendix~\ref{app:normal-trace}.
\end{proof}

\paragraph{The residual estimate.}\label{subsec:residual}

\begin{proposition}\label{prop:residual}
Let $g^\lambda$ be as in \eqref{eq:g-def}. There are $C_T,c>0$ such that for every $\tau\in(0,T)$,
\begin{equation}\label{eq:residual-L2}
  \|g^\lambda\|_{L^2(\tau,T;L^2(\Gamma))}
  \ \le\ C_T\,(1+\tau^{-1/2})
  \bigl(1+\|\hh_0\|_\lambda\bigr)\,
  \bigl(\lambda^{1/4}+e^{-c/\lambda}\bigr),
\end{equation}
and, if in addition $E(0)<\infty$ (finite initial weighted gradient energy; cf.\ Section~\ref{subsec:lifted}),
\begin{equation}\label{eq:residual-L2-zero}
  \|g^\lambda\|_{L^2(0,T;L^2(\Gamma))}
  \ \le\ C_T\,\bigl(1+\|\hh_0\|_\lambda+E(0)^{1/2}\bigr)\,
  \bigl(\lambda^{1/4}+e^{-c/\lambda}\bigr).
\end{equation}
\end{proposition}

\begin{proof}
Fix $\theta\in C^\infty([0,\delta'])$ with $\theta(0)=1$, $\theta\equiv0$ on $[2\delta'/3,\delta']$, $\theta'$ supported in $[\delta'/3,2\delta'/3]$. Fix a time $t$ for which Lemma~\ref{lem:transmission} holds. By Lemma~\ref{lem:transmission}, $s\mapsto\theta(s)F(t,\cdot,s)$ is in $W^{1,1}((0,\delta');L^2(\Gamma))$ with value $g^\lambda(t,\cdot)$ at $s=0$ and $0$ at $s=\delta'$. Hence, in $L^2(\Gamma)$,
\[
  g^\lambda(t,\cdot)
  =-\int_0^{\delta'}\partial_s\bigl(\theta F\bigr)\dd s
  =-\int_0^{\delta'}\theta'F\dd s
  -\int_0^{\delta'}\theta\,\partial_sF\dd s .
\]
Substituting \eqref{eq:dsF} in the $L^2(\Gamma)$-valued integral gives
\begin{equation*}
  g^\lambda(t,y)
  =-\int_0^{\delta'}\theta\,a_\lambda\,\partial_tq\dd s
  -\int_0^{\delta'}\theta'\,a_\lambda\,(\partial_sq+\beta_sq)\dd s
  +\int_0^{\delta'}\theta\,\divG\!\left[
  a_\lambda\bigl(\mathsf B\nG q+\beta_\Gamma q\bigr)\right]\dd s .
\end{equation*}
Every integral is over the exterior side, where
\[
\int_0^{\delta'}a_\lambda(y,s)\dd s\le
C\int_0^\infty e^{-s^2/2\lambda}\dd s\le C\sqrt\lambda .
\]
We estimate the three terms in space--time.

\emph{Term 1.} By Cauchy--Schwarz in $s$ with weight $a_\lambda$,
\[
  \Bigl|\int_0^{\delta'}\theta a_\lambda\,\partial_tq\dd s\Bigr|^2
  \le\Bigl(\int_0^{\delta'}a_\lambda\dd s\Bigr)
  \int_0^{\delta'}a_\lambda\,|\partial_tq|^2\dd s
  \le C\sqrt\lambda\int_0^{\delta'}a_\lambda\,|\partial_tq|^2\dd s .
\]
% Integrating over $\Gamma$ and $t\in(\tau,T)$,
Using \(dx=m(y,s)\,d\sigma(y)\,ds\), we have
\[
\int_\Gamma\int_0^{\delta'}
a_\lambda(y,s)|\partial_tq(t,y,s)|^2\,ds\,d\sigma(y)
\le \|\partial_t\widehat h(t)\|_\lambda^2.
\]
Integrating in time and taking square roots, we obtain
\[
  \bigl\|\text{Term 1}\bigr\|_{L^2((\tau,T)\times\Gamma)}
  \le C\lambda^{1/4}
  \Bigl(\int_\tau^T\|\partial_t\hh\|_\lambda^2\Bigr)^{1/2}
  \le C\lambda^{1/4}(1+\tau^{-1/2})
  \bigl(1+\|\hh_0\|_\lambda\bigr),
\]
by Lemma~\ref{lem:smoothing}(i).

\emph{Term 2.} On the support of $\theta'$, $s\ge\delta'/3$, so $a_\lambda\le Ce^{-\delta'^2/36\lambda}\,a_\lambda^{1/2}$ pointwise (write $a_\lambda=a_\lambda^{1/2}a_\lambda^{1/2}$ and bound one factor by $e^{-s^2/4\lambda}\le e^{-\delta'^2/36\lambda}$). Hence
\[
  \bigl\|\text{Term 2}\bigr\|_{L^2((\tau,T)\times\Gamma)}
  \le Ce^{-c/\lambda}
  \Bigl(\int_\tau^T\!\!\bigl(E+\|\hh\|_\lambda^2\bigr)\Bigr)^{1/2}
  \le Ce^{-c/\lambda}\bigl(1+\|\hh_0\|_\lambda\bigr) .
\]

\emph{Term 3.} By \eqref{eq:tangential-div-bound} and the same weighted Cauchy--Schwarz inequality in $s$ as for \emph{Term 1},
\[
  \bigl\|\text{Term 3}(t)\bigr\|_{L^2(\Gamma)}
  \le C\lambda^{1/4}
  \Bigl(\bigl\|\nG^2q\bigr\|_{L^2(a_\lambda)}
  +\bigl\|\nG q\bigr\|_{L^2(a_\lambda)}
  +\|q\|_{L^2(a_\lambda)}\Bigr)(t),
\]
where the norms are taken over the exterior part $\Gamma\times(0,\delta')$ (bounded by the full-collar norms). Integrating in $t\in(\tau,T)$ and invoking Lemma~\ref{lem:tangential}, Lemma~\ref{lem:smoothing}(i), and Lemma~\ref{lem:basic-energy} bounds this contribution by $C\lambda^{1/4}(1+\tau^{-1/2})(1+\|\hh_0\|_\lambda)$. Summing the three estimates proves \eqref{eq:residual-L2}; using Lemma~\ref{lem:smoothing}(ii) instead of (i) gives \eqref{eq:residual-L2-zero}.
\end{proof}

\subsection{Auxiliary initial density and stability}\label{subsec:lifted}

The hard datum $\hh_0=\rho_0\one_K$ has a jump across $\Gamma$ whenever $\rho_0|_{\Gamma}\not\equiv0$, so $E(0)=\infty$ and \eqref{eq:residual-L2-zero} is not directly available. We therefore compare it with an auxiliary solution whose initial relative density has finite weighted gradient energy. The resulting $L^2(\muh)$ initial mismatch is of order $\lambda^{1/4}$.

Let $\rho_0\in H^1(K)$ with $0\le\rho_0\le M$. By Sobolev extension and truncation there is $H_0\in H^1(\Rd)\cap L^\infty(\Rd)$ with $H_0=\rho_0$ on $K$, $0\le H_0\le M$, and $\|H_0\|_{H^1(\Rd)}\le C_K\|\rho_0\|_{H^1(K)}\le C_KN_0$. Define the normalized auxiliary relative density and the corresponding solution
\[
  \bar h_0:=\kappa_\lambda H_0,\qquad
  \kappa_\lambda^{-1}:=\int_{\Rd}H_0\dd\muh,\qquad
  \bar\rho^\lambda:=\text{penalized solution with datum }\bar h_0\,\wl,
\]
and $\bar h:=\bar\rho^\lambda/\wl$. Since $\int H_0\dd\muh=\int_K\rho_0+\int_{K^c}H_0\dd\muh =1+O(M\sqrt\lambda)$ by Lemma~\ref{lem:layer}, we have $|\kappa_\lambda-1|\le CM\sqrt\lambda$, and
\begin{equation*}
  \|\hh_0-\bar h_0\|_\lambda^2
  =|1-\kappa_\lambda|^2\int_K\rho_0^2\dd x
  +\kappa_\lambda^2\int_{K^c}H_0^2\dd\muh
  \le CN_0^2\lambda+CM^2\,\muh(K^c)
  \le CN_0^2\sqrt\lambda,
\end{equation*}
so
\begin{equation}\label{eq:init-mismatch-final}
  \|\hh_0-\bar h_0\|_\lambda\ \le\ C\,N_0\,\lambda^{1/4},
  \qquad
  \int_{\Rd}|\nabla\bar h_0|^2\dd\muh
  \le C\|H_0\|_{H^1(\Rd)}^2\le CN_0^2 .
\end{equation}
Moreover, $0<\kappa_\lambda\le1$, so $\|\bar h_0\|_\lambda^2\le M\int_{\Rd}\bar h_0\dd\muh=M$. The auxiliary solution has drift $b[\bar\rho^\lambda(t)]\in\mathcal B$ and $\bar h_0\wl\le CM\wl$; hence Proposition~\ref{thm:barrier} gives the same $O(\sqrt\lambda)$ exterior mass bound.

\begin{lemma}
\label{lem:pen-stability}
Let $\rho^{(1)},\rho^{(2)}$ be penalized solutions with data $\rho^{(i)}_0=h^{(i)}_0\wl$, $h^{(i)}_0\in L^2(\muh)$, $\rho^{(i)}_0\in\mcP(\Rd)$. Then, with $z:=h^{(1)}-h^{(2)}$,
\[
  \sup_{t\le T}\|z_t\|_\lambda^2
  +\int_0^T\!\!\int|\nabla z|^2\dd\muh\dd t
  \ \le\ C\,\|z_0\|_\lambda^2,
\]
where $C$ depends on $T,B,L_W,K$, and $\sup_{t\le T}\|h^{(2)}_t\|_\lambda$, but not on $\lambda$.
\end{lemma}

\begin{proof}
Subtracting the two ratio equations \eqref{eq:ratio-eq},
\[
  \wl\partial_tz=\divg\Bigl(\wl\nabla z+\wl z\,b[\rho^{(1)}_t]
  +\wl h^{(2)}\bigl(b[\rho^{(1)}_t]-b[\rho^{(2)}_t]\bigr)\Bigr).
\]
Testing with $z$ gives
\[
  \frac12\frac{\dd}{\dd t}\|z\|_\lambda^2+\|\nabla z\|_{L^2(\muh)}^2
  \le B\|z\|_\lambda\|\nabla z\|_{L^2(\muh)}
  +\|b[\rho^{(1)}_t]-b[\rho^{(2)}_t]\|_\infty
  \|h^{(2)}\|_\lambda\|\nabla z\|_{L^2(\muh)},
\]
and, by \eqref{eq:drift-bounds} and Cauchy--Schwarz, $\|b[\rho^{(1)}_t]-b[\rho^{(2)}_t]\|_\infty \le L_W\|\rho^{(1)}_t-\rho^{(2)}_t\|_{L^1} =L_W\int|z|\dd\muh\le L_W\muh(\Rd)^{1/2}\|z\|_\lambda\le C\|z\|_\lambda$. Young's inequality followed by Gr\"onwall's inequality yields the assertion.
\end{proof}

Combining \eqref{eq:init-mismatch-final} with Lemma~\ref{lem:pen-stability}, with $h^{(1)}=\hh$ and $h^{(2)}=\bar h$, gives
\begin{equation}\label{eq:hard-vs-lifted}
  \sup_{t\le T}\bigl\|\rho^\lambda(t)-\bar\rho^\lambda(t)\bigr\|_{L^2(K)}
  +\bigl\|\rho^\lambda-\bar\rho^\lambda\bigr\|_{L^2(0,T;H^1(K))}
  \ \le\ C\,N_0\,\lambda^{1/4},
\end{equation}
since on $K$ the weight is $1$. Moreover, by \eqref{eq:residual-L2-zero} applied to $\bar\rho^\lambda$ (whose datum has $\int_{\Rd}|\nabla\bar h_0|^2\dd\muh\le CN_0^2$),
\begin{equation}\label{eq:lifted-residual}
  \bigl\|\bar g^\lambda\bigr\|_{L^2(0,T;L^2(\Gamma))}
  \le CN_0\bigl(\lambda^{1/4}+e^{-c/\lambda}\bigr),
\end{equation}
where $\bar g^\lambda$ is the trace of $\bar J^\lambda_K:=\nabla\bar\rho^\lambda+\bar\rho^\lambda b[\bar\rho^\lambda]$ provided by Lemma~\ref{lem:transmission}.

\subsection{Comparison on the domain and completion of the proof}
\label{subsec:finite-comparison}

\begin{proposition}
\label{prop:error-estimate}
Let $\rho$ be the weak solution of \eqref{eq:ref-pde}. Let $\tilde\rho^\lambda$ be a penalized solution with boundary residual $\tilde g^\lambda\in L^2(0,T;L^2(\Gamma))$, and suppose
\[
  \sup_{t\le T}\|\tilde\rho^\lambda(t)\|_{L^2(K)}\le C_0,
  \qquad
  \sup_{t\le T}\int_{K^c}\tilde\rho^\lambda(t)\,\dd x\le C_0\sqrt\lambda.
\]
Set $u:=\tilde\rho^\lambda|_K-\rho$. Then
\begin{equation}\label{eq:error-estimate}
  \sup_{t\le T}\|u(t)\|_{L^2(K)}^2
  +\int_0^T\|\nabla u\|_{L^2(K)}^2\dd t
  \ \le\
  C_T(C_0)\Bigl(\|u_0\|_{L^2(K)}^2
  +\lambda
  +\|\tilde g^\lambda\|_{L^2((0,T)\times\Gamma)}^2\Bigr).
\end{equation}
\end{proposition}

\begin{proof}
For a.e.\ $t$, the Green identity \eqref{eq:trace-green} states precisely that, for $\varphi\in H^1(K)$,
\[
  \langle\partial_t\tilde\rho^\lambda,\varphi\rangle
  =-\int_K\bigl(\nabla\tilde\rho^\lambda
  +\tilde\rho^\lambda b[\tilde\rho^\lambda]\bigr)
  \cdot\nabla\varphi\dd x+\int_\Gamma\tilde g^\lambda\varphi\dd\Hm .
\]
Thus $\tilde\rho^\lambda|_K$ satisfies the unpenalized equation with normal flux $\tilde g^\lambda$. Subtracting \eqref{eq:weak-form},
\[
  \langle\partial_tu,\varphi\rangle
  +\int_K\Bigl(\nabla u+u\,b[\rho(t)]
  +\tilde\rho^\lambda\bigl(b[\tilde\rho^\lambda(t)]-b[\rho(t)]\bigr)\Bigr)
  \cdot\nabla\varphi\dd x
  =\int_\Gamma\tilde g^\lambda\varphi\dd\Hm.
\]
The interior fluxes are in \(L^2((0,T)\times K)\), and the trace theorem controls the boundary term by \(C\|\tilde g^\lambda(t)\|_{L^2(\Gamma)}\|\varphi\|_{H^1(K)}\). Hence \(\partial_tu\in L^2(0,T;H^1(K)')\). Together with \(u\in L^2(0,T;H^1(K))\), this justifies the test \(\varphi=u(t)\) by the Lions--Magenes chain rule and gives
\[
\begin{aligned}
  \frac12\frac{\dd}{\dd t}\|u\|_{L^2(K)}^2
  +\|\nabla u\|_{L^2(K)}^2
  &\le B\|u\|_{L^2(K)}\|\nabla u\|_{L^2(K)}\\
  &\quad+\bigl\|\tilde\rho^\lambda
  (b[\tilde\rho^\lambda(t)]-b[\rho(t)])\bigr\|_{L^2(K)}
  \|\nabla u\|_{L^2(K)}
  +\Bigl|\int_\Gamma\tilde g^\lambda u\dd\Hm\Bigr| .
\end{aligned}
\]
By \eqref{eq:drift-bounds} and the assumed exterior-mass bound,
\[
  \bigl\|b[\tilde\rho^\lambda(t)]-b[\rho(t)]\bigr\|_\infty
  \le L_W\Bigl(\|u(t)\|_{L^1(K)}
  +\int_{K^c}\tilde\rho^\lambda(t)\Bigr)
  \le C\bigl(\|u(t)\|_{L^2(K)}+\sqrt\lambda\bigr),
\]
and the assumed bound $\|\tilde\rho^\lambda(t)\|_{L^2(K)}\le C_0$ gives $\|\tilde\rho^\lambda (b[\tilde\rho^\lambda(t)]-b[\rho(t)])\|_{L^2(K)} \le C(\|u\|_{L^2(K)}+\sqrt\lambda)$. The boundary term is estimated by the multiplicative trace inequality and Young's inequality:
\[
  \Bigl|\int_\Gamma\tilde g^\lambda u\dd\Hm\Bigr|
  \le\|\tilde g^\lambda\|_{L^2(\Gamma)}\|u\|_{L^2(\Gamma)}
  \le\frac14\|\nabla u\|_{L^2(K)}^2
  +C\|u\|_{L^2(K)}^2+C\|\tilde g^\lambda\|_{L^2(\Gamma)}^2 .
\]
Absorbing the gradient terms,
\[
  \frac{\dd}{\dd t}\|u\|_{L^2(K)}^2+\|\nabla u\|_{L^2(K)}^2
  \le C\|u\|_{L^2(K)}^2+C\lambda
  +C\|\tilde g^\lambda(t)\|_{L^2(\Gamma)}^2.
\]
An application of Gr\"onwall's inequality gives \eqref{eq:error-estimate}.
\end{proof}

\begin{proof}[Proof of Theorem~\ref{thm:strong-main}]
\emph{Hard initialization.} Let $\rho^\lambda$ denote the solution with the hard datum $\rho_0\one_K$. Since its residual estimate is available only away from $t=0$, we first use the auxiliary flow $\bar\rho^\lambda$ constructed in Section~\ref{subsec:lifted}. By Lemma~\ref{lem:interior} (with \(\|\bar h_0\|_\lambda^2\le M\)) and Proposition~\ref{thm:barrier}, this flow satisfies the two bounds assumed in Proposition~\ref{prop:error-estimate}, with \(C_0\) independent of \(\lambda\). Its initial mismatch with the no-flux datum on $K$ is $|1-\kappa_\lambda|\,\|\rho_0\|_{L^2(K)}\le CN_0\sqrt\lambda$, and its boundary residual satisfies \eqref{eq:lifted-residual}. Proposition \ref{prop:error-estimate} therefore gives
\[
  \sup_{t\le T}\|\bar\rho^\lambda(t)-\rho(t)\|_{L^2(K)}
  +\|\bar\rho^\lambda-\rho\|_{L^2(0,T;H^1(K))}
  \le CN_0\bigl(\lambda^{1/4}+e^{-c/\lambda}\bigr)
  \le CN_0\lambda^{1/4}.
\]
The stability comparison \eqref{eq:hard-vs-lifted} and the triangle inequality now yield the first two terms of \eqref{eq:strong-rate} for the hard solution $\rho^\lambda$.

Set $u=\rho^\lambda|_K-\rho$. For the corresponding interior fluxes,
\[
  J^\lambda_K-J=\nabla u+u\,b[\rho(t)]
  +\rho^\lambda\bigl(b[\rho^\lambda(t)]-b[\rho(t)]\bigr),
\]
and the three terms are bounded in $L^2((0,T)\times K)$ by, respectively, the $H^1$ part of \eqref{eq:strong-rate}, $B$ times its $L^2$ part, and $\sup_t\|\rho^\lambda(t)\|_{L^2(K)}\cdot \|b[\rho^\lambda]-b[\rho]\|_{L^2(0,T;L^\infty)} \le C_T(\sup_t\|u(t)\|_{L^2(K)}+\sqrt\lambda) \le C_T\lambda^{1/4}$. This proves the flux term in \eqref{eq:strong-rate}. The boundary estimate \eqref{eq:bdry-rate} follows by applying \eqref{eq:residual-L2} directly to the hard solution and absorbing \(e^{-c/\lambda}\) into \(C_T\lambda^{1/4}\).

\emph{Auxiliary initial density.} For the auxiliary flow, the only change relevant to the boundary estimate is that \eqref{eq:init-mismatch-final} gives $\int_{\Rd}|\nabla\bar h_0|^2\dd\muh\le CN_0^2$. Lemma~\ref{lem:smoothing}(ii), and hence the second part of Proposition~\ref{prop:residual}, is therefore valid from $t=0$. Estimate \eqref{eq:lifted-residual}, with the exponential term absorbed into \(C_T\lambda^{1/4}\), proves the final assertion of the theorem.
\end{proof}

% {\color{red}Estimate \eqref{eq:lifted-residual}, with the exponential term absorbed
% into \(C_T\lambda^{1/4}\), gives the boundary estimate from \(t=0\)
% for the auxiliary flow.
% }

\begin{remark}
The interaction enters through the time regularity of the drift, the stability estimate for the penalized flow, and the drift mismatch in Proposition~\ref{prop:error-estimate}. These terms do not alter the power \(\lambda^{1/4}\), but the constants resulting from Gr\"onwall's inequality depend on \(T\). Section~\ref{sec:counterexample} shows that this dependence cannot be removed for general interactions.
\end{remark}

\subsection{Sharpness under hard initialization}\label{subsec:sharpness}

The rate $\lambda^{1/4}$ in \eqref{eq:strong-rate} is optimal even for the linear heat equation. The obstruction is the initial layer created by extending a nonzero boundary value by zero.

\begin{proposition}[Sharpness under hard initialization]\label{prop:sharpness}
Let $K=[0,1]$, $\Phi=W=0$, and $\rho_0\equiv1$ on $K$, extended by zero to $\R$. Let $\rho^\lambda$ and $\rho$ be the corresponding penalized and no-flux solutions. Then $\rho\equiv1$, $J=0$, and $J_K^\lambda=\partial_x\rho^\lambda$ on $K$. There exist constants $c_1,c_2>0$ such that, for every $T>0$ and all sufficiently small $\lambda>0$,
\begin{align*}
 \|\rho^\lambda-\rho\|_{L^\infty(0,T;L^2(K))}
 &\ge c_1\lambda^{1/4},\\
 \|\partial_x(\rho^\lambda-\rho)\|_{L^2((0,T)\times K)}
 =\|J_K^\lambda-J\|_{L^2((0,T)\times K)}
 &\ge c_2\lambda^{1/4}.
\end{align*}
\end{proposition}

\begin{proof}
Put $L=\lambda^{-1/2}$, $v(\tau,y)=\rho^\lambda(\lambda\tau,1+\sqrt\lambda\,y)$, and $d_L(y)=y-\proj_{[-L,0]}y$. Then
\[
 v_\tau=v_{yy}+\partial_y(d_Lv),\qquad v(0)=\one_{[-L,0]}.
\]
The profile $w_L=e^{-\dist(\cdot,[-L,0])^2/2}$ is stationary, since $w_L'+d_Lw_L=0$. As $\one_{[-L,0]}\le w_L$, comparison gives $0\le v\le w_L$, hence $\|d_Lv\|_\infty\le e^{-1/2}$ and $\|d_Lv\|_1\le2$. Let $G_\tau(y)=(4\pi\tau)^{-1/2}e^{-y^2/(4\tau)}$ and $H_L(\tau)=G_\tau*\one_{[-L,0]}$. Duhamel's formula gives
\[
 v(\tau)=H_L(\tau)
 +\int_0^\tau\partial_yG_{\tau-s}*(d_Lv(s))\,\dd s.
\]
The boundedness and integrability of $d_Lv$ justify this identity in the distributional sense. Since $\|\partial_yG_s\|_1=(\pi s)^{-1/2}$, uniformly in $L$,
\begin{equation}\label{eq:sharpness-heat}
 \|v(\tau)-H_L(\tau)\|_\infty
 \le e^{-1/2}\int_0^\tau\|\partial_yG_s\|_1\,\dd s
 =2e^{-1/2}\sqrt{\tau/\pi}\le\sqrt\tau.
\end{equation}
Fix $\tau_0>0$ so small that $2\tau_0+\sqrt{\tau_0}\le1/8$, and choose $0<a<1$ with $a/\sqrt{2\pi\tau_0}\le1/8$. For $L\ge2$ and $\tau\in I:=[\tau_0/2,\tau_0]$,
\[
 H_L(\tau,y)\le\frac12+\frac{|y|}{\sqrt{4\pi\tau}}
 \quad(-a\le y\le0),\qquad
 H_L(\tau,-1)\ge1-2\tau.
\]
The second inequality follows from the Gaussian second moment, since $[-2,0]\subset[-L,0]$. Together with \eqref{eq:sharpness-heat}, these estimates imply
\[
 v(\tau,y)\le\frac34\quad(-a\le y\le0),\qquad
 v(\tau,-1)\ge\frac78.
\]
Consequently, for almost every $\tau\in I$,
\[
 \int_{-a}^{0}|v(\tau,y)-1|^2\,\dd y\ge\frac a{16},
 \qquad
 \int_{-1}^{0}|v_y(\tau,y)|^2\,\dd y
 \ge|v(\tau,-1)-v(\tau,0)|^2\ge\frac1{64},
\]
where endpoint values are Sobolev traces and the second estimate is the spatial Cauchy--Schwarz inequality. For $\lambda$ small enough that $L\ge2$ and $\lambda\tau_0\le T$, rescaling gives
\[
 \mathop{\rm ess\,sup}_{0<t<T}
 \|\rho^\lambda(t)-1\|_{L^2(K)}^2
 \ge\frac a{16}\sqrt\lambda,\qquad
 \int_0^T\!\int_K|\partial_x\rho^\lambda|^2\,\dd x\,\dd t
 \ge\frac{\tau_0}{128}\sqrt\lambda.
\]
Taking $c_1=\sqrt a/4$ and $c_2=\sqrt{\tau_0/128}$ proves the claim, since $\partial_x\rho=J=0$ and $J_K^\lambda=\partial_x\rho^\lambda$ on $K$.
\end{proof}

The lower bounds arise at times $t\asymp\lambda$ from the jump in the zero-extended initial density. They establish optimality for all three error terms in \eqref{eq:strong-rate} on time intervals starting at $t=0$. On fixed intervals $[\tau,T]$ with $\tau>0$, the initial layer is excluded, and better rates for these errors and the normal-flux trace may be possible.

\section{Proof of Theorem~\ref{thm:linear-global}}
\label{sec:linear-global}

We now assume \(W\equiv0\). To prove Theorem~\ref{thm:linear-global}, we use the reversible structure both in the whole-space estimates and in the comparison on \(K\). We write
\begin{equation*}
 \begin{aligned}
  Z_0&:=\int_K e^{-\Phi(x)}\,\dd x,
  &\pi(x)&:=Z_0^{-1}e^{-\Phi(x)} \quad (x\in K),\\
  Z_\lambda&:=\int_{\Rd}e^{-\Phi(x)-U_\lambda(x)}\,\dd x,
  &\pi_\lambda(x)&:=Z_\lambda^{-1}e^{-\Phi(x)-U_\lambda(x)}
  \quad (x\in\Rd).
 \end{aligned}
\end{equation*}
The corresponding interior and hard fluxes are \(J_K^\lambda(t)=\nabla\rho^\lambda(t)+\rho^\lambda(t)\nabla\Phi\) and \(J(t)=\nabla\rho(t)+\rho(t)\nabla\Phi\).

Here the natural relative density is \(\rho^\lambda/\pi_\lambda\): division by the Gibbs density absorbs all drift terms, and equilibrium corresponds to the constant \(1\). Dissipation is used in two distinct steps. On \(\Rd\), the \(\lambda\)-uniform Poincar\'e inequality and the energy identities for the auxiliary flow control the deviation from equilibrium and the derivatives needed in the boundary-layer argument over \((0,\infty)\); see \eqref{eq:linear-global-energy}. The exterior representation from Section~\ref{sec:strong} then yields the global residual bound of Lemma~\ref{lem:global-linear-residual}. Subtracting the equilibrium is essential for these time-integrable estimates.

On \(K\), a second use of dissipation prevents this small boundary forcing from accumulating in the density error. We divide this error by \(\pi\); the resulting variable has \(\pi\)-mean equal to minus the exterior mass, which is \(O(\sqrt\lambda)\) uniformly in time, so the Poincar\'e inequality for \(\pi\) on \(K\) applies to it up to an additive error \(O(\lambda)\). This yields the dissipative comparison of Lemma~\ref{lem:linear-dissipative}, with constants independent of time. Thus the proof combines global control of the residual with a dissipative, rather than merely finite-time, comparison.

\subsection{Uniform Gibbs estimates and auxiliary initial data}\label{subsec:linear-lifted}

\begin{lemma}\label{lem:uniform-gibbs-gap}
There exist \(C,C_P>0\) and \(\lambda_0>0\) such that, for \(0<\lambda\le\lambda_0\),
\begin{equation}\label{eq:gibbs-layer}
  0\le\int_{K^c}\pi_\lambda\,\dd x
  =1-\frac{Z_0}{Z_\lambda}\le C\sqrt\lambda,
  \qquad \pi_\lambda|_K=\frac{Z_0}{Z_\lambda}\pi .
\end{equation}
Moreover, for every \(v\in H^1_{\mathrm{loc}}(\Rd)\) satisfying \(\int_{\Rd}(v^2+|\nabla v|^2)\pi_\lambda\,\dd x<\infty\),
\begin{equation}\label{eq:uniform-poincare}
  \int_{\Rd}\left|v-\int_{\Rd}v\pi_\lambda\,\dd x\right|^2
  \pi_\lambda\,\dd x
  \le C_P\int_{\Rd}|\nabla v|^2\pi_\lambda\,\dd x .
\end{equation}
\end{lemma}

\begin{proof}
Since \(\Phi\) is Lipschitz and \(\proj_Kx\in K\), \(e^{-\Phi(x)}\le C e^{\|\nabla\Phi\|_\infty\dist(x,K)}\). Lemma~\ref{lem:layer} therefore gives
\[
  0\le Z_\lambda-Z_0
  =\int_{K^c}e^{-\Phi(x)-\dist(x,K)^2/(2\lambda)}\,\dd x
  \le C\sqrt\lambda.
\]
Dividing by \(Z_\lambda\ge Z_0>0\) proves the exterior-mass estimate in \eqref{eq:gibbs-layer}; the identity on \(K\) follows directly from \(U_\lambda=0\).

We next prove a uniform estimate along the exterior normal rays. By convexity, every \(x\in K^c\) can be written uniquely as \(x=\Psi(y,s):=y+s\vec n(y)\), with \(y=\proj_Kx\in\Gamma\) and \(s=\dist(x,K)>0\). If \(\tau_j\) is a principal direction at \(y\) with curvature \(\kappa_j(y)\), then
\[
  D\Psi(y,s)\tau_j=(1+s\kappa_j(y))\tau_j,
  \qquad \partial_s\Psi(y,s)=\vec n(y),
  \qquad m(y,s)=\prod_{j=1}^{d-1}(1+s\kappa_j(y)).
\]
Here \(m\) is the Jacobian of \(\Psi\). The nonnegative curvatures are uniformly bounded, and hence the following estimates hold uniformly in the displayed variables:
\[
  1\le m(y,s)\le C(1+s)^{d-1},
  \qquad \frac{m(y,r+u)}{m(y,r)}\le C(1+u)^{d-1}.
\]

For fixed \(y\in\Gamma\), define the one-dimensional weight \(\omega_{\lambda,y}(s):= e^{-\Phi(\Psi(y,s))-s^2/(2\lambda)}m(y,s)\). We claim that, uniformly in \(y\), \(r\ge0\), and sufficiently small \(\lambda\),
\begin{equation}\label{eq:normal-ray-bounds}
  \int_0^\infty\omega_{\lambda,y}(s)\,\dd s
  \le C\sqrt\lambda\,e^{-\Phi(y)},
  \qquad
  \int_r^\infty s\omega_{\lambda,y}(s)\,\dd s
  \le C\lambda\,\omega_{\lambda,y}(r).
\end{equation}
The first estimate follows from the Lipschitz bound for \(\Phi\), the polynomial growth of \(m\), and the substitution \(s=\sqrt\lambda z\). For the second, write \(s=r+u\). The Jacobian-ratio estimate and the same Lipschitz bound give
\begin{align*}
 \frac{1}{\omega_{\lambda,y}(r)}
 \int_r^\infty s\omega_{\lambda,y}(s)\,\dd s
 &\le C\int_0^\infty(r+u)(1+u)^{d-1}
 e^{-ru/\lambda-u^2/(2\lambda)+\|\nabla\Phi\|_\infty u}\,\dd u\\
 &=C\lambda\int_0^\infty
 \left(\frac r{\sqrt\lambda}+z\right)(1+\sqrt\lambda z)^{d-1}
 e^{-(r/\sqrt\lambda)z-z^2/2+\|\nabla\Phi\|_\infty\sqrt\lambda z}
 \,\dd z\le C\lambda.
\end{align*}
The last bound follows from \(\sup_{a\ge0}\int_0^\infty(a+z)e^{-az-z^2/4}\,\dd z<\infty\).

For smooth \(v\), the fundamental theorem of calculus and Cauchy--Schwarz give, along each normal ray,
\[
  |v(\Psi(y,s))|^2\le2|\Tr v(y)|^2
  +2s\int_0^s|\nabla v(\Psi(y,r))|^2\,\dd r.
\]
Multiplying by \(\omega_{\lambda,y}(s)\), integrating in \(s\), and using Tonelli's theorem and \eqref{eq:normal-ray-bounds}, we obtain
\begin{equation}\label{eq:normal-hardy}
  \int_0^\infty|v(\Psi(y,s))|^2\omega_{\lambda,y}(s)\,\dd s
  \le C\sqrt\lambda\,e^{-\Phi(y)}|\Tr v(y)|^2
  +C\lambda\int_0^\infty|\nabla v(\Psi(y,s))|^2
  \omega_{\lambda,y}(s)\,\dd s.
\end{equation}

We now combine this estimate with the spectral gap on \(K\). Subtracting the \(\pi_\lambda\)-mean, we may assume \(\int_{\Rd}v\pi_\lambda\,\dd x=0\). Since \(\pi_\lambda=(Z_0/Z_\lambda)\pi\) on \(K\), Cauchy--Schwarz and \eqref{eq:gibbs-layer} imply
\begin{equation*}
  \left|\int_Kv\pi\,\dd x\right|^2
  \le C\sqrt\lambda\int_{K^c}v^2\pi_\lambda\,\dd x.
\end{equation*}
After integrating \eqref{eq:normal-hardy} over \(\Gamma\) and dividing by \(Z_\lambda\), the trace theorem and the Poincar\'e inequality for \(\pi\) give
\[
 \begin{aligned}
  \int_{K^c}v^2\pi_\lambda\,\dd x
  &\le C\sqrt\lambda\left(
      \int_K|\nabla v|^2\pi\,\dd x
      +\left|\int_Kv\pi\,\dd x\right|^2\right)
      +C\lambda\int_{K^c}|\nabla v|^2\pi_\lambda\,\dd x\\
  &\le C\sqrt\lambda\int_{\Rd}|\nabla v|^2\pi_\lambda\,\dd x
      +C\lambda\int_{K^c}v^2\pi_\lambda\,\dd x.
 \end{aligned}
\]
For small \(\lambda\), the last term is absorbed. Thus
\[
  \int_{K^c}v^2\pi_\lambda\,\dd x
  \le C\sqrt\lambda\int_{\Rd}|\nabla v|^2\pi_\lambda\,\dd x,
  \qquad
  \left|\int_Kv\pi\,\dd x\right|^2
  \le C\lambda\int_{\Rd}|\nabla v|^2\pi_\lambda\,\dd x.
\]
Together with the Poincar\'e inequality on \(K\), these two estimates give \eqref{eq:uniform-poincare} for smooth functions. Truncation and density give the stated result.
\end{proof}

\paragraph{Auxiliary initial density.}
The global dissipative estimate is formulated for densities relative to the invariant measure \(\pi_\lambda\). We therefore construct an auxiliary initial density by extension and weighting with the invariant density, adapting the construction of Section~\ref{subsec:lifted}. Set \(f_0:=\rho_0/\pi\) on \(K\). Since \(\pi^{-1}\in W^{1,\infty}(K)\), the assumptions on \(\rho_0\) imply \(f_0\in H^1(K)\cap L^\infty(K)\). The Sobolev extension theorem, followed by truncation and multiplication by a cutoff, provides a nonnegative, compactly supported \(F_0\in H^1(\Rd)\cap L^\infty(\Rd)\) such that \(F_0|_K=f_0\). Define
\begin{equation}\label{eq:linear-lift}
  c_\lambda^{-1}:=\int_{\Rd}F_0\pi_\lambda\,\dd x,
  \qquad \bar f_0^\lambda:=c_\lambda F_0,
  \qquad \bar\rho_0^\lambda:=\pi_\lambda\bar f_0^\lambda.
\end{equation}
Because \(\int_Kf_0\pi\,\dd x=1\),
\[
  c_\lambda^{-1}=\frac{Z_0}{Z_\lambda}
  +\int_{K^c}F_0\pi_\lambda\,\dd x,
  \qquad 0\le\int_{K^c}F_0\pi_\lambda\,\dd x\le C\sqrt\lambda.
\]
Together with \eqref{eq:gibbs-layer}, this gives
\begin{equation}\label{eq:linear-lift-bounds}
  |c_\lambda-1|\le C\sqrt\lambda,\qquad
  \|\bar f_0^\lambda\|_\infty+
  \int_{\Rd}|\nabla\bar f_0^\lambda|^2\pi_\lambda\,\dd x\le C.
\end{equation}
On \(K\), the auxiliary density differs from \(\rho_0\) only by the constant factor \(c_\lambda Z_0/Z_\lambda\). Its normalization and the preceding exterior-layer estimate show that this factor differs from one by at most \(C\sqrt\lambda\). Hence \(\|\bar\rho_0^\lambda|_K-\rho_0\|_{L^2(K)}\le C\sqrt\lambda\).

For hard initialization, the relative initial density is \((Z_\lambda/Z_0)f_0\one_K\). Its squared distance from the auxiliary relative density separates into an interior normalization error and an exterior layer error:
\begin{align}
 \left\|\frac{Z_\lambda}{Z_0}f_0\one_K
       -\bar f_0^\lambda\right\|_{L^2(\pi_\lambda)}^2
 &=\frac{Z_0}{Z_\lambda}
   \left|\frac{Z_\lambda}{Z_0}-c_\lambda\right|^2
   \|f_0\|_{L^2(\pi)}^2
   +c_\lambda^2\int_{K^c}F_0^2\pi_\lambda\,\dd x \notag\\
 &\le C\lambda+C\sqrt\lambda\le C\sqrt\lambda.
 \label{eq:linear-lift-mismatch}
\end{align}
Thus the corresponding \(L^2(\pi_\lambda)\)-norm is \(O(\lambda^{1/4})\); the slower rate comes from the exterior layer.

\subsection{Global boundary residual}

\begin{lemma}\label{lem:global-linear-residual}
Let \(\bar\rho^\lambda\) be the penalized solution of the linear equation with initial datum \eqref{eq:linear-lift}, and let \(\bar f^\lambda:=\bar\rho^\lambda/\pi_\lambda\). Its interior normal-flux trace \(\bar g^\lambda:=(\nabla\bar\rho^\lambda +\bar\rho^\lambda\nabla\Phi)\cdot\vec n\) satisfies
\begin{align}
 \|\bar g^\lambda\|_{L^2(0,\infty;L^2(\Gamma))}
 &\le C\lambda^{1/4},\label{eq:global-linear-residual}\\
 \sup_{t\ge0}\int_{K^c}\bar\rho^\lambda(t)\,\dd x
 &\le C\sqrt\lambda.\label{eq:global-linear-exterior}
\end{align}
\end{lemma}

\begin{proof}
Set \(v^\lambda(t):=\bar f^\lambda(t)-1\). Subtracting the equilibrium is needed for a time-integrable $L^2$ estimate; it leaves the flux unchanged. Since \(\bar\rho^\lambda(t)\) has unit mass, \(v^\lambda(t)\) has zero \(\pi_\lambda\)-mean for every \(t\ge0\), and it solves
\begin{equation}\label{eq:linear-relative}
  \pi_\lambda\partial_tv^\lambda
  =\divg(\pi_\lambda\nabla v^\lambda).
\end{equation}
The initial bound \eqref{eq:linear-lift-bounds} gives \(\|v^\lambda(0)\|_{L^2(\pi_\lambda)}^2 +\|\nabla v^\lambda(0)\|_{L^2(\pi_\lambda)}^2\le C\). Testing \eqref{eq:linear-relative} successively with \(v^\lambda(t)\) and \(\partial_tv^\lambda(t)\) yields
\begin{align*}
 \frac12\frac{\dd}{\dd t}\|v^\lambda(t)\|_{L^2(\pi_\lambda)}^2
 +\|\nabla v^\lambda(t)\|_{L^2(\pi_\lambda)}^2&=0,\\
 \frac12\frac{\dd}{\dd t}\|\nabla v^\lambda(t)\|_{L^2(\pi_\lambda)}^2
 +\|\partial_tv^\lambda(t)\|_{L^2(\pi_\lambda)}^2&=0.
\end{align*}
The second identity follows by integrating by parts in space. Both identities are justified for the auxiliary relative density by the standard self-adjoint semigroup approximation for form-domain data. After integration in time, the uniform Poincar\'e inequality \eqref{eq:uniform-poincare} gives
\begin{equation}\label{eq:linear-global-energy}
  \int_0^\infty
  \left(
  \|v^\lambda(t)\|_{L^2(\pi_\lambda)}^2
  +\|\nabla v^\lambda(t)\|_{L^2(\pi_\lambda)}^2
  +\|\partial_tv^\lambda(t)\|_{L^2(\pi_\lambda)}^2
  \right)\dd t
  +\sup_{t\ge0}\|\nabla v^\lambda(t)\|_{L^2(\pi_\lambda)}^2
  \le C.
\end{equation}

On \(K\), the identity \(\nabla\pi_\lambda+\pi_\lambda\nabla\Phi=0\) gives
\[
 \bar g^\lambda
 =\pi_\lambda\nabla v^\lambda\cdot\vec n
 \qquad\text{on }\Gamma.
\]
Write \(\widehat v^\lambda(t,y,s):=v^\lambda(t,\Psi(y,s))\) and set
\[
 \begin{aligned}
 \widetilde a_\lambda(y,s)
 &:=\pi_\lambda(\Psi(y,s))m(y,s)
 =Z_\lambda^{-1}e^{-\Phi(\Psi(y,s))}a_\lambda(y,s),\\
 \widetilde F&:=\widetilde a_\lambda\,\partial_s\widehat v^\lambda.
 \end{aligned}
\]
The collar equation is now
\[
 \widetilde a_\lambda\,\partial_t\widehat v^\lambda
 =\partial_s\widetilde F+
 \divG\bigl(\widetilde a_\lambda\mathsf B\nG\widehat v^\lambda\bigr),
\]
and \(\widetilde F(t,y,0)=\bar g^\lambda(t,y)\) by the same trace identification as in Lemma~\ref{lem:transmission}. On the fixed collar, the factor \(Z_\lambda^{-1}e^{-\Phi\circ\Psi}\) is bounded above and below by positive constants and has bounded first derivatives, uniformly in \(\lambda\). Hence
\[
 \sup_{y\in\Gamma}\int_0^\delta\widetilde a_\lambda(y,s)\,\dd s
 \le C\sqrt\lambda,\qquad
 \|\nG\log\widetilde a_\lambda\|_\infty\le C.
\]
We can therefore repeat the collar representation and the proof of Proposition~\ref{prop:residual} with this weight. The coefficients are time-independent, so the constants in this spatial argument do not grow with the time horizon. We use the global energy bound \eqref{eq:linear-global-energy}, rather than pass to infinite time in an estimate with constant $C_T$. The difference-quotient argument of Lemma~\ref{lem:tangential} and the residual representation give
\begin{align*}
 \|\bar g^\lambda\|_{L^2(0,\infty;L^2(\Gamma))}
 \le{}&C\lambda^{1/4}\bigl(
 \|\partial_tv^\lambda\|_{L^2(0,\infty;L^2(\pi_\lambda))}
 +\|v^\lambda\|_{L^2(0,\infty;L^2(\pi_\lambda))}
 +\|\nabla v^\lambda\|_{L^2(0,\infty;L^2(\pi_\lambda))}\bigr)\\
 &+Ce^{-c/\lambda}
 \|\nabla v^\lambda\|_{L^2(0,\infty;L^2(\pi_\lambda))}.
\end{align*}
Here the tangential second derivatives are controlled by the same right-hand side, exactly as in Lemma~\ref{lem:tangential}. The global energy estimate \eqref{eq:linear-global-energy} now proves \eqref{eq:global-linear-residual}.

Finally, the maximum principle gives \(0\le\bar f^\lambda(t)\le\|\bar f_0^\lambda\|_\infty\) for every \(t\ge0\). Hence \eqref{eq:linear-lift-bounds} and \eqref{eq:gibbs-layer} imply
\[
  \sup_{t\ge0}\int_{K^c}\bar\rho^\lambda(t)\,\dd x
  \le \|\bar f_0^\lambda\|_\infty
      \int_{K^c}\pi_\lambda\,\dd x
  \le C\sqrt\lambda,
\]
which proves \eqref{eq:global-linear-exterior}.
\end{proof}

\subsection{Dissipative comparison and completion of the proof}

\begin{lemma}\label{lem:linear-dissipative}
Let \(\widetilde\rho^\lambda\) be a penalized solution of the linear equation such that \(\sup_{t\ge0}\int_{K^c}\widetilde\rho^\lambda(t)\,\dd x \le C_0\sqrt\lambda\), and suppose that its interior normal-flux trace is \(g^\lambda\in L^2(0,\infty;L^2(\Gamma))\). Set
\[
 u(t):=\widetilde\rho^\lambda(t)|_K-\rho(t),\qquad
 \widetilde J_K^\lambda(t):=
 \nabla\widetilde\rho^\lambda(t)+\widetilde\rho^\lambda(t)\nabla\Phi.
\]
Then
\begin{align}
 &\sup_{t\ge0}\|u(t)\|_{L^2(K)}^2
 +\sup_{s\ge0}\int_s^{s+1}
   \left(\|u(t)\|_{H^1(K)}^2+
   \|\widetilde J_K^\lambda(t)-J(t)\|_{L^2(K)}^2\right)\dd t
 \notag\\
 &\qquad\le
 C\left(\|u(0)\|_{L^2(K)}^2+\lambda+
 \|g^\lambda\|_{L^2(0,\infty;L^2(\Gamma))}^2\right).
 \label{eq:linear-dissipative}
\end{align}
\end{lemma}

\begin{proof}
Let \(r(t):=u(t)/\pi\). Since \(\nabla u+u\nabla\Phi=\pi\nabla r\), subtracting the two equations and testing the resulting error equation with \(r(t)\) gives
\begin{equation}\label{eq:linear-error-energy}
  \frac12\frac{\dd}{\dd t}\int_Kr(t)^2\pi\,\dd x
  +\int_K|\nabla r(t)|^2\pi\,\dd x
  =\int_\Gamma g^\lambda(t)r(t)\,\dd\Hm .
\end{equation}
Mass conservation yields
\[
  \int_Kr(t)\pi\,\dd x=\int_Ku(t)\,\dd x
  =-\int_{K^c}\widetilde\rho^\lambda(t)\,\dd x,
\]
whose absolute value is at most \(C_0\sqrt\lambda\). The weighted Poincar\'e inequality on \(K\), followed by the trace theorem, therefore gives
\begin{equation}\label{eq:linear-error-poincare}
  \int_Kr(t)^2\pi\,\dd x+\|r(t)\|_{L^2(\Gamma)}^2
  \le C\int_K|\nabla r(t)|^2\pi\,\dd x+C\lambda.
\end{equation}
Using \eqref{eq:linear-error-poincare} and Young's inequality to estimate the right-hand side of \eqref{eq:linear-error-energy}, and using part of the gradient term to control the weighted \(L^2\)-term once more through \eqref{eq:linear-error-poincare}, we obtain constants \(c,\kappa>0\) such that
\begin{align*}
 &\frac{\dd}{\dd t}\int_Kr(t)^2\pi\,\dd x
 +\kappa\int_Kr(t)^2\pi\,\dd x
 +c\int_K|\nabla r(t)|^2\pi\,\dd x\\
 &\qquad\le
 C\lambda+C\|g^\lambda(t)\|_{L^2(\Gamma)}^2.
\end{align*}
Variation of constants now gives
\[
  \sup_{t\ge0}\int_Kr(t)^2\pi\,\dd x
  \le C\left(
  \int_Kr(0)^2\pi\,\dd x+\lambda
  +\|g^\lambda\|_{L^2(0,\infty;L^2(\Gamma))}^2
  \right).
\]
Integrating the same differential inequality over \([s,s+1]\) and using this supremum bound yields
\[
  \sup_{s\ge0}\int_s^{s+1}\int_K
  |\nabla r(t)|^2\pi\,\dd x\dd t
  \le C\left(
  \int_Kr(0)^2\pi\,\dd x+\lambda
  +\|g^\lambda\|_{L^2(0,\infty;L^2(\Gamma))}^2
  \right).
\]
Finally,
\[
  \widetilde J_K^\lambda(t)-J(t)=\pi\nabla r(t),\qquad
  \nabla u(t)=\pi\nabla r(t)-u(t)\nabla\Phi,
\]
and \(\pi\) is bounded above and below by positive constants on \(K\). The preceding estimates therefore imply \eqref{eq:linear-dissipative}.
\end{proof}

\begin{proof}[Proof of Theorem~\ref{thm:linear-global}]
Applying Lemma~\ref{lem:linear-dissipative} to the auxiliary solution and using \eqref{eq:linear-lift-bounds}, \eqref{eq:global-linear-residual}, and \eqref{eq:global-linear-exterior} gives, with \(\bar J_K^\lambda(t):=\nabla\bar\rho^\lambda(t) +\bar\rho^\lambda(t)\nabla\Phi\),
\begin{align}
 &\sup_{t\ge0}\|\bar\rho^\lambda(t)|_K-\rho(t)\|_{L^2(K)}^2 \notag\\
 &\quad+\sup_{s\ge0}\int_s^{s+1}\left(
 \|\bar\rho^\lambda(t)|_K-\rho(t)\|_{H^1(K)}^2
 +\|\bar J_K^\lambda(t)-J(t)\|_{L^2(K)}^2\right)\dd t
 \le C\sqrt\lambda .
 \label{eq:linear-lifted-conclusion}
\end{align}

It remains to compare the auxiliary solution with the solution having hard initialization. Let \(f_h^\lambda(t):=\rho^\lambda(t)/\pi_\lambda\). The relative densities \(f_h^\lambda\) and \(\bar f^\lambda\) evolve under the same reversible semigroup. Hence \eqref{eq:linear-lift-mismatch}, the energy identity, and \eqref{eq:uniform-poincare} imply
\[
  \sup_{t\ge0}\|f_h^\lambda(t)-\bar f^\lambda(t)\|_{L^2(\pi_\lambda)}^2
  +\int_0^\infty\left(
   \|f_h^\lambda(t)-\bar f^\lambda(t)\|_{L^2(\pi_\lambda)}^2
   +\|\nabla(f_h^\lambda(t)-\bar f^\lambda(t))\|_{L^2(\pi_\lambda)}^2
   \right)\dd t
  \le C\sqrt\lambda .
\]
On \(K\),
\[
  \rho^\lambda-\bar\rho^\lambda
  =\pi_\lambda(f_h^\lambda-\bar f^\lambda),\qquad
  J_K^\lambda-\bar J_K^\lambda
  =\pi_\lambda\nabla(f_h^\lambda-\bar f^\lambda).
\]
Since \(\pi_\lambda\) and its first derivatives are uniformly bounded on \(K\), the same squared estimate holds for the interior \(H^1\) and flux differences. Combining it with \eqref{eq:linear-lifted-conclusion} and taking square roots proves \eqref{eq:linear-global-rate}.

The maximum principle gives \(0\le f_h^\lambda(t)\le\|f_h^\lambda(0)\|_\infty\le C\), and therefore
\[
  \sup_{t\ge0}\int_{K^c}\rho^\lambda(t)\,\dd x\le C\sqrt\lambda.
\]
Decomposing the full-space $L^1$ error over $K$ and $K^c$ gives
\[
\begin{aligned}
  \|\rho^\lambda(t)-\rho(t)\one_K\|_{L^1(\Rd)}
  &=\|\rho^\lambda(t)|_K-\rho(t)\|_{L^1(K)}
     +\int_{K^c}\rho^\lambda(t)\,\dd x\\
  &\le |K|^{1/2}\|\rho^\lambda(t)|_K-\rho(t)\|_{L^2(K)}
     +\int_{K^c}\rho^\lambda(t)\,\dd x.
\end{aligned}
\]
Combining \eqref{eq:linear-global-rate} with the preceding exterior-mass bound yields \(C\lambda^{1/4}+C\sqrt\lambda\), uniformly in time. Since \(\sqrt\lambda\le\lambda^{1/4}\) for \(0<\lambda\le1\), this proves \eqref{eq:linear-global-L1}.
\end{proof}

\section{Proof of Theorem~\ref{thm:counterexample}}\label{sec:counterexample}

We prove Theorem~\ref{thm:counterexample} by constructing an unstable equilibrium for which the finite-time approximation remains valid but uniform-in-time convergence fails.

\subsection{The model and its unstable equilibrium}

Let
\[
  K=[-\pi,\pi],\qquad
  U_\lambda(x)=\frac{(|x|-\pi)_+^2}{2\lambda},
  \qquad \lambda>0,
\]
and fix
\begin{equation*}
  \mathcal W(z)=-a\cos z,
  \qquad a>2.
\end{equation*}
The penalized equation on $\R$ is
\begin{equation}\label{eq:soft-full}
  \partial_t\rho^\lambda
  =\partial_x\!\left(
      \partial_x\rho^\lambda
      +U_\lambda'\rho^\lambda
      +\rho^\lambda\,\partial_x(\mathcal W*\rho^\lambda(t))
    \right),
\end{equation}
with the hard initialization
\begin{equation*}
  \rho_0^\lambda(x)=\bar\rho\,\one_K(x),
  \qquad
  \bar\rho:=\frac{1}{2\pi}.
\end{equation*}
The no-flux equation is
\begin{equation}\label{eq:hard-full}
\begin{cases}
  \displaystyle
  \partial_t\rho
  =\partial_x\!\left(
      \partial_x\rho
      +\rho\,\partial_x(\mathcal W*\rho(t))
    \right),&x\in(-\pi,\pi),\\[0.4em]
  \displaystyle
  \left(
      \partial_x\rho
      +\rho\,\partial_x(\mathcal W*\rho(t))
  \right)\big|_{x=\pm\pi}=0,
\end{cases}
\end{equation}
with initial datum $\rho_0=\bar\rho$. In \eqref{eq:hard-full}, the convolution is taken against the density on $K$.

We use the penalized and weak solution of Definitions~\ref{def:pen-sol} and \ref{def:weak-sol}; both solutions conserve mass. For every fixed $\lambda>0$, uniqueness of the penalized problem and invariance under $x\mapsto-x$ imply
\begin{equation*}
  \rho^\lambda(t)(-x)=\rho^\lambda(t,x),
  \qquad t\geq0.
\end{equation*}
Define
\begin{equation*}
  m_\lambda(t):=\int_\R\cos x\,\rho^\lambda(t,x)\dd x.
\end{equation*}
Even symmetry and the addition formula for sine and cosine give $\partial_x(\mathcal W*\rho^\lambda(t))=a m_\lambda(t)\sin x$. Consequently, \eqref{eq:soft-full} reduces to
\begin{equation}\label{eq:soft-reduced}
  \partial_t\rho^\lambda
  =\partial_x\!\left(
      \partial_x\rho^\lambda
      +U_\lambda'\rho^\lambda
      +a m_\lambda(t)\sin x\,\rho^\lambda
    \right).
\end{equation}
For the narrow-convergence conclusion, extend the observable \(\sin^2x\) from \(K\) to \(\R\) by
\begin{equation}\label{eq:observable-extension}
  \widetilde\varphi(x)=
  \begin{cases}
    \sin^2x,&|x|\leq\pi,\\
    0,&|x|>\pi.
  \end{cases}
\end{equation}
Since both the value and the first derivative of $\sin^2x$ vanish at $x=\pm\pi$, one has $\widetilde\varphi\in C_b^1(\R)$. Moreover, \(\int_K\sin^2x\,\bar\rho\dd x=1/2\).

\paragraph{The hard equilibrium.}

The interaction drift of the uniform density vanishes. Hence $\bar\rho=1/(2\pi)$ is stationary for \eqref{eq:hard-full}, and uniqueness gives
\begin{equation}\label{eq:hard-constant}
  \rho(t)\equiv\bar\rho.
\end{equation}

On even mass-zero perturbations, the linearized operator is
\begin{equation*}
  Lu=u''+a\bar\rho
  \left(\int_K\cos y\,u(y)\dd y\right)\cos x,
\end{equation*}
with $u'(\pm\pi)=0$. In particular,
\begin{equation}\label{eq:unstable-eigenvalue}
  L(\cos x)=\left(\frac a2-1\right)\cos x.
\end{equation}
Thus the hard stationary density is linearly unstable whenever $a>2$.

\subsection{The cosine moment}

Define
\begin{equation*}
  p_\lambda(t)=\int_\R\sin x\,U_\lambda'(x)
  \rho^\lambda(t,x)\dd x,\qquad
  \sigma_\lambda(t)
  =a\int_\R\sin^2x\,\rho^\lambda(t,x)\dd x-1.
\end{equation*}
Testing \eqref{eq:soft-reduced} against $\cos x$ gives
\begin{equation}\label{eq:exact-moment}
  m_\lambda'(t)=\sigma_\lambda(t)m_\lambda(t)+p_\lambda(t).
\end{equation}
The identity follows rigorously by testing with cutoffs $\chi_R\cos x$ and letting $R\to\infty$, using Lemma~\ref{lem:basic-energy} and \eqref{eq:domination}.

Since \(\rho^\lambda(t)\) is a probability density,
\begin{equation}\label{eq:sigma-uniform-bound}
  -1\leq\sigma_\lambda(t)\leq a-1,
  \qquad \abs{\sigma_\lambda(t)}\leq a-1.
\end{equation}
Even symmetry gives the exterior representation
\begin{equation}\label{eq:p-exterior}
  p_\lambda(t)
  =-\frac{2}{\lambda}
  \int_0^\infty s\sin s\,
  \rho^\lambda(t,\pi+s)\dd s.
\end{equation}
The contribution of the first exterior Gaussian layer, $0<s<\pi$, to the cosine moment is therefore negative. The regions where $\sin s<0$ contribute with the opposite sign, but they begin at distance $\pi$ and are exponentially small on the logarithmic time scale considered below.

\subsection{Boundary-layer forcing and the initial unstable amplitude}

For every fixed \(T\), Proposition~\ref{thm:barrier}, applied to the relative density \(\hh^\lambda=\rho^\lambda e^{U_\lambda}\), gives
\begin{equation}\label{eq:counter-ratio-bound}
  0\le \hh^\lambda(t,x)\le C_Te^{\beta\dist(x,K)}
  \qquad (0\le t\le T).
\end{equation}
Here \(\beta>0\) is the constant in Proposition~\ref{thm:barrier}; it depends only on \(K\), the dimension, and the drift class, and is independent of \(T\) and \(\lambda\). We determine the leading boundary-layer forcing.

\begin{lemma}\label{lem:boundary-forcing}
For every $0<\tau<T<\infty$,
\begin{equation*}
  \frac{p_\lambda}{\sqrt\lambda}
  \longrightarrow
  -\frac{1}{\sqrt{2\pi}}
  \qquad\text{in }L^1(\tau,T).
\end{equation*}
\end{lemma}

\begin{proof}
Set
\[
  h_+^\lambda(t,s):=\hh^\lambda(t,\pi+s),
  \qquad s\geq0.
\]
We first identify its value at the boundary. For each fixed \(\lambda>0\), the weight \(e^{-U_\lambda}\) and its reciprocal are bounded on every bounded neighborhood of \(x=\pi\). The weighted energy regularity of \(\hh^\lambda(t)\) therefore implies that \(\hh^\lambda(t)\in H^1(\pi-\delta,\pi+\delta)\) for every \(\delta>0\) and for a.e.\ \(t>0\). In one dimension this function has a continuous representative near \(\pi\), so its left and right traces coincide. Since \(U_\lambda(\pi)=0\),
\[
  \hh^\lambda(t,\pi^+)=\hh^\lambda(t,\pi^-)=\rho^\lambda(t,\pi^-).
\]
By \eqref{eq:strong-rate} and the one-dimensional trace theorem,
\begin{equation}\label{eq:trace-convergence}
  h_+^\lambda(\cdot,0)
  =\rho^\lambda(\cdot,\pi^-)
  \longrightarrow\bar\rho
  \qquad\text{in }L^2(\tau,T).
\end{equation}
For fixed $z>0$, weighted Cauchy--Schwarz gives, for almost every $t$,
\begin{equation*}
 \begin{aligned}
  \abs{h_+^\lambda(t,\sqrt\lambda z)-h_+^\lambda(t,0)}^2
  &\leq
  \left(
    \int_0^{\sqrt\lambda z}
    e^{-s^2/(2\lambda)}
    \abs{\partial_sh_+^\lambda(t,s)}^2\dd s
  \right)
  \left(
    \int_0^{\sqrt\lambda z}e^{s^2/(2\lambda)}\dd s
  \right)
  \\
  &=
  \sqrt\lambda
  \left(\int_0^z e^{r^2/2}\dd r\right)
  \int_0^{\sqrt\lambda z}
    e^{-s^2/(2\lambda)}
    \abs{\partial_sh_+^\lambda(t,s)}^2\dd s.
 \end{aligned}
\end{equation*}
After integration over $t\in(\tau,T)$, Lemma~\ref{lem:basic-energy} implies
\begin{equation}\label{eq:normal-comparison-L2}
  \norm{
    h_+^\lambda(\cdot,\sqrt\lambda z)
    -h_+^\lambda(\cdot,0)
  }_{L^2(\tau,T)}^2
  \leq C_{T,z}\sqrt\lambda.
\end{equation}
Combining \eqref{eq:trace-convergence} and \eqref{eq:normal-comparison-L2},
\begin{equation}\label{eq:profile-convergence}
  h_+^\lambda(\cdot,\sqrt\lambda z)
  \longrightarrow\bar\rho
  \qquad\text{in }L^2(\tau,T)
\end{equation}
for every fixed $z>0$.

Using $\rho^\lambda=e^{-U_\lambda}\hh^\lambda$ in \eqref{eq:p-exterior} and substituting $s=\sqrt\lambda z$, we obtain
\begin{equation}\label{eq:p-scaled}
  \frac{p_\lambda(t)}{\sqrt\lambda}
  =-2\int_0^\infty
  z\,\frac{\sin(\sqrt\lambda z)}{\sqrt\lambda}
  e^{-z^2/2}
  h_+^\lambda(t,\sqrt\lambda z)\dd z.
\end{equation}
For $\lambda\leq1$,
\[
  z\abs{\frac{\sin(\sqrt\lambda z)}{\sqrt\lambda}}
  \leq z^2.
\]
On the fixed interval $t\in[\tau,T]$, the ratio bound \eqref{eq:counter-ratio-bound} gives
\[
  h_+^\lambda(t,\sqrt\lambda z)
  \leq C_Te^{\beta\sqrt\lambda z}
  \leq C_Te^{\beta z}.
\]
Thus the integrand in \eqref{eq:p-scaled}, viewed in $L^1(\tau,T)$, is dominated by an integrable multiple of
\[
  z^2e^{-z^2/2+\beta z}.
\]
Using \eqref{eq:profile-convergence}, dominated convergence in the $z$-integral yields
\begin{align*}
  \frac{p_\lambda}{\sqrt\lambda}
  &\longrightarrow
  -2\bar\rho\int_0^\infty z^2e^{-z^2/2}\dd z=-\frac{1}{\pi}\sqrt{\frac{\pi}{2}}
  =-\frac{1}{\sqrt{2\pi}}
\end{align*}
in $L^1(\tau,T)$.
\end{proof}

The following lemma controls the oscillatory tail that can make \(p_\lambda\) positive.

\begin{lemma}
\label{lem:positive-tail}
Let \(p_\lambda^+(t):=\max\{p_\lambda(t),0\}\). There exist \(c,C,\lambda_0>0\) such that, for \(0<\lambda\le\lambda_0\),
\begin{equation*}
  p_\lambda^+(t)
  \leq Ce^{\gamma t-c/\lambda},
  \qquad t\geq0.
\end{equation*}
Here \(\gamma\) is the constant in Proposition~\ref{thm:barrier}.
\end{lemma}

\begin{proof}
Positive contributions in \eqref{eq:p-exterior} require $s\geq\pi$. Using \eqref{eq:domination}, for sufficiently small $\lambda$,
\begin{align*}
  p_\lambda^+(t)
  &\leq\frac{2}{\lambda}\int_\pi^\infty
    s\rho^\lambda(t,\pi+s)\dd s
  \leq\frac{C}{\lambda}e^{\gamma t}
  \int_\pi^\infty s e^{-s^2/(4\lambda)}\dd s\leq Ce^{\gamma t-\pi^2/(4\lambda)}.
\end{align*}
\end{proof}

Consequently, for every fixed \(\ell,N>0\), \(\sup_{0\le t\le\ell|\log\lambda|}p_\lambda^+(t)=O(\lambda^N)\). Thus the positive part is negligible on every logarithmic time interval.

\begin{corollary}
\label{cor:negative-moment}
There exist $c_0>0$ and $\lambda_0>0$ such that
\begin{equation*}
  m_\lambda(2)\leq-c_0\sqrt\lambda
  \qquad\text{for }0<\lambda<\lambda_0.
\end{equation*}
\end{corollary}

\begin{proof}
Write $p_\lambda=p_\lambda^+-p_\lambda^-$. Lemmas \ref{lem:boundary-forcing} and \ref{lem:positive-tail} imply
\[
  \int_1^2p_\lambda^-(t)\dd t\geq c\sqrt\lambda,
  \qquad \int_0^2p_\lambda^+(t)\dd t\leq Ce^{-c/\lambda}.
\]
Since $m_\lambda(0)=0$, variation of constants in \eqref{eq:exact-moment} yields
\[
  m_\lambda(2)
  =\int_0^2
  \exp\!\left(\int_r^2\sigma_\lambda(\xi)\dd\xi\right)
  p_\lambda(r)\dd r.
\]
Together with \eqref{eq:sigma-uniform-bound}, this gives $m_\lambda(2)\le-e^{-(a-1)}\int_1^2p_\lambda^-+ e^{2(a-1)}\int_0^2p_\lambda^+\le-c_0\sqrt\lambda$.
\end{proof}

\subsection{Weak separation on the logarithmic time scale}

The moment equation supplies the growing mode, but its growth coefficient still depends on the solution. The observable $\sin^2x$ is chosen because it controls this coefficient from below. If all its unit-time averages were small, the negative cosine moment would grow beyond the bound $|m_\lambda|\le1$. The argument uses the sign of the exterior contribution, not a small exterior-mass estimate on the entire logarithmic interval.

Define the interior observable
\begin{equation*}
  \mathcal O_\lambda(t)
  :=\int_K\sin^2x\,\bigl(\rho^\lambda(t,x)-\bar\rho\bigr)\dd x.
\end{equation*}
Since the exterior contribution to the integral defining \(\sigma_\lambda\) is nonnegative,
\begin{equation}\label{eq:sigma-lower-pointwise}
  \sigma_\lambda(t)
  \geq a\int_K\sin^2x\,\rho^\lambda(t,x)\dd x-1
  =\mu_a+a\mathcal O_\lambda(t),
\end{equation}
where
\begin{equation*}
  \mu_a:=\frac a2-1>0.
\end{equation*}
Choose \(\eta_a:=\mu_a/(4a)\), and set \(\kappa_a:=\mu_a-a\eta_a=3\mu_a/4\). If, for some \(T>3\),
\begin{equation*}
  \sup_{2\leq s\leq T-1}
  \abs{\int_s^{s+1}\mathcal O_\lambda(t)\dd t}
  \leq\eta_a,
\end{equation*}
then partitioning $[r,t]$ into unit intervals and using \(\mathcal O_\lambda\ge-1/2\) on the remainder gives, for \(2\le r\le t\le T\),
\[
  \int_r^t\mathcal O_\lambda(\xi)\dd\xi
  \geq-\eta_a(t-r)-\frac12.
\]
Indeed, every complete unit interval contributes at least $-\eta_a$, while the remaining interval has length less than one and contributes at least $-1/2$. Combining this estimate with \eqref{eq:sigma-lower-pointwise} yields
\begin{equation}\label{eq:sigma-average-lower}
  \int_r^t\sigma_\lambda(\xi)\dd\xi
  \geq\kappa_a(t-r)-\frac a2.
\end{equation}

\begin{proposition}
\label{prop:log-loss}
There is $C_a>0$ such that, for every sufficiently small $\lambda$, one can find $s_\lambda$ satisfying \eqref{eq:s-lambda-bound} and
\begin{equation}\label{eq:localized-window-failure}
  \abs{\int_{s_\lambda}^{s_\lambda+1}\mathcal O_\lambda(t)\dd t}
  \geq\eta_a.
\end{equation}
\end{proposition}

\begin{proof}
Let $c_0$ be the constant in Corollary~\ref{cor:negative-moment} and define
\begin{equation*}
  T_\lambda
  :=2+\frac{1}{\kappa_a}
  \left(
    \frac12|\log\lambda|
    +\frac a2
    +\log\frac{4}{c_0}
  \right).
\end{equation*}
For small $\lambda$, $T_\lambda>3$ and $T_\lambda\leq C_a(1+|\log\lambda|)$.

Assume, for contradiction, that
\begin{equation}\label{eq:no-window-before-T}
  \sup_{2\leq s\leq T_\lambda-1}
  \abs{\int_s^{s+1}\mathcal O_\lambda(t)\dd t}
  <\eta_a.
\end{equation}
Variation of constants in \eqref{eq:exact-moment}, starting from time $2$, gives
\begin{equation*}
  m_\lambda(t)
  =\exp\!\left(\int_2^t\sigma_\lambda(\xi)\dd\xi\right)
    m_\lambda(2)
  +\int_2^t
    \exp\!\left(\int_r^t\sigma_\lambda(\xi)\dd\xi\right)
    p_\lambda(r)\dd r.
\end{equation*}
The negative part of $p_\lambda$ only decreases the right-hand side. Using Corollary~\ref{cor:negative-moment}, \eqref{eq:sigma-average-lower}, the upper bound \(\sigma_\lambda\leq a-1\), and Lemma~\ref{lem:positive-tail}, we obtain, for $2\leq t\leq T_\lambda$,
\begin{align}
  m_\lambda(t)
  &\leq
  -c_0\sqrt\lambda
  \exp\!\left(\kappa_a(t-2)-\frac a2\right)
  +\int_2^t e^{(a-1)(t-r)}p_\lambda^+(r)\dd r
  \notag\\
  &\leq
  -c_0\sqrt\lambda
  \exp\!\left(\kappa_a(t-2)-\frac a2\right)
  +C(1+t)e^{(a-1+\gamma)t-c/\lambda}.
  \label{eq:m-upper-log}
\end{align}
At $t=T_\lambda$, the first term on the right-hand side is exactly $-4$. Since $T_\lambda=O(|\log\lambda|)$,
\begin{equation*}
  (1+T_\lambda)
  e^{(a-1+\gamma)T_\lambda-c/\lambda}
  \longrightarrow0.
\end{equation*}
Thus \eqref{eq:m-upper-log} gives $m_\lambda(T_\lambda)<-2$ for all sufficiently small $\lambda$. This is impossible because $\rho^\lambda(t)$ is a probability density and $\abs{\cos x}\leq1$, so
\begin{equation*}
  \abs{m_\lambda(t)}\leq1
  \qquad\text{for every }t\geq0.
\end{equation*}
The window integral is continuous in its starting time, since $\rho^\lambda\in C([0,T_\lambda];L^1(\R))$. Its absolute value attains a maximum on $[2,T_\lambda-1]$. Therefore \eqref{eq:no-window-before-T} is false, and a time $s_\lambda\in[2,T_\lambda-1]$ satisfying \eqref{eq:localized-window-failure} exists.
\end{proof}

\begin{proof}[Proof of Theorem~\ref{thm:counterexample}]
The identity $\rho(t)\equiv\bar\rho$ follows from \eqref{eq:hard-constant}. Proposition~\ref{prop:log-loss} gives \eqref{eq:s-lambda-bound} and \eqref{eq:window-weak-main}. Since the time interval in \eqref{eq:window-weak-main} has length one, for every sufficiently small $\lambda$,
\[
  \sup_{t\ge0}
  \left|\int_K \sin^2x\,
  \bigl(\rho^\lambda(t,x)-\bar\rho\bigr)\dd x\right|
  \ge \eta_a.
\]
Thus the pairing with the fixed function $\sin^2x\in L^2(K)$ does not converge to zero uniformly in time. This rules out uniform-in-time weak convergence in $L^2(K)$, and hence also uniform-in-time strong convergence. For narrow convergence, consider $\widetilde\varphi\in C_b(\R)$ which is obtained by extending $\sin^2 x$ as $0$ outside of $K$. Since
\[
  \int_{\R}\widetilde\varphi(x)
  \bigl(\rho^\lambda(t,x)-\bar\rho\one_K(x)\bigr)\dd x
  =
  \int_K\sin^2x\,
  \bigl(\rho^\lambda(t,x)-\bar\rho\bigr)\dd x,
\]
this rules out uniform-in-time narrow convergence to $\bar\rho\one_K(x)\dd x$. Finally $s_\lambda \to \infty$ follows by Theorem~\ref{thm:strong-main} as stated.
\end{proof}

\appendix

% ==================================================
\section{Smoothing estimate by Galerkin method}
\label{app:smoothing}
% ==================================================

\begin{proof}[Detailed proof of Lemma~\ref{lem:smoothing}]
Fix \(0<\lambda\le1\) and freeze \(b^\lambda\) along the already constructed solution \(\rho^\lambda\). Choose nested finite-dimensional spaces \(V_n\subset H^1(\muh)\) whose union is dense in \(H^1(\muh)\). For an initial approximation \(\hh_{n,0}\in V_n\), let \(\hh_n\in C^1([0,T];V_n)\) solve
\begin{equation}\label{eq:smoothing-galerkin}
 (\partial_t\hh_n,v)_\lambda+(\nabla\hh_n,\nabla v)_\lambda
 +(\hh_nb^\lambda,\nabla v)_\lambda=0
 \quad\text{for every }v\in V_n,
\end{equation}
where \((\cdot,\cdot)_\lambda\) denotes integration against \(\muh\). In a basis of \(V_n\), this is a linear nonautonomous system of ordinary differential equations. For part~(ii), take \(\hh_{n,0}\to\hh_0\) in \(H^1(\muh)\); convergence in \(L^2(\muh)\) suffices for part~(i).

Set
\[
 D_n(t):=\|\partial_t\hh_n(t)\|_\lambda^2,\qquad
 E_n(t):=\|\nabla\hh_n(t)\|_\lambda^2.
\]
The basic estimate gives, uniformly in \(n\) and \(\lambda\),
\begin{equation}\label{eq:smoothing-basic-galerkin}
 \sup_{t\le T}\|\hh_n(t)\|_\lambda^2+\int_0^TE_n(t)\,\dd t\le C.
\end{equation}
Let \(G(t):=\sup_x|\dot b^\lambda(t,x)|\). By Lemma~\ref{lem:drift-time-reg}, \(G\) is uniformly bounded in \(L^2(0,T)\). For basis functions \(e_i,e_j\in V_n\), the map \(t\mapsto\int_{\Rd}e_i b^\lambda(t)\cdot\nabla e_j\,\dd\muh\) is absolutely continuous, and its derivative is obtained by replacing \(b^\lambda\) with \(\dot b^\lambda\). Hence \(Y_n(t):=\int_{\Rd}\hh_n b^\lambda\cdot\nabla\hh_n\,\dd\muh\) belongs to \(W^{1,1}(0,T)\).

Testing \eqref{eq:smoothing-galerkin} by \(\partial_t\hh_n\) and using the product rule for \(Y_n\) yields
\begin{equation}\label{eq:smoothing-differential}
 \frac12D_n+\frac12E_n'+Y_n'
 \le\frac{B^2}{2}E_n
 +G(t)\|\hh_n(t)\|_\lambda E_n(t)^{1/2}.
\end{equation}
Moreover,
\begin{equation}\label{eq:smoothing-Y-coercivity}
 |Y_n(t)|\le B\|\hh_n(t)\|_\lambda E_n(t)^{1/2}
 \le\frac14E_n(t)+B^2\|\hh_n(t)\|_\lambda^2.
\end{equation}

For \(H^1(\muh)\) data, integrate \eqref{eq:smoothing-differential} on \((0,t)\), use \eqref{eq:smoothing-Y-coercivity} at the endpoints, and estimate
\[
 \int_0^tG\|\hh_n\|_\lambda E_n^{1/2}\,\dd s
 \le\sup_{s\le T}\|\hh_n(s)\|_\lambda\|G\|_{L^2(0,T)}
 \left(\int_0^TE_n\,\dd s\right)^{1/2}.
\]
Together with \eqref{eq:smoothing-basic-galerkin}, this gives \(\int_0^tD_n(s)\,\dd s+E_n(t)\le C(1+E_n(0))\). Passing to the limit proves part~(ii).

For \(L^2(\muh)\) data, multiply \eqref{eq:smoothing-differential} by \(t\) and integrate. The initial energy then disappears, and integration by parts gives
\begin{align*}
 \frac12\int_0^t sD_n(s)\,\dd s+\frac t2E_n(t)+tY_n(t)
 \le{}&\frac12\int_0^tE_n(s)\,\dd s+\int_0^tY_n(s)\,\dd s\\
 &+\frac{B^2}{2}\int_0^t sE_n(s)\,\dd s
 +\int_0^t sG(s)\|\hh_n(s)\|_\lambda E_n(s)^{1/2}\,\dd s .
\end{align*}
The terms containing \(Y_n\) are controlled by \eqref{eq:smoothing-Y-coercivity}; the last term is bounded by
\[
 T\sup_{s\le T}\|\hh_n(s)\|_\lambda\|G\|_{L^2(0,T)}
 \left(\int_0^TE_n(s)\,\dd s\right)^{1/2}.
\]
It follows that
\begin{equation}\label{eq:weighted-smoothing-final}
 \int_0^T sD_n(s)\,\dd s
 +\operatorname*{ess\,sup}_{0<t\le T}tE_n(t)\le C.
\end{equation}

The preceding bounds yield, along a subsequence,
\[
 \hh_n\rightharpoonup\widetilde h
 \ \text{in }L^2(0,T;H^1(\muh)),\qquad
 \hh_n\stackrel{*}{\rightharpoonup}\widetilde h
 \ \text{in }L^\infty(0,T;L^2(\muh)),
\]
and, for every \(\tau>0\), \(\partial_t\hh_n\rightharpoonup\partial_t\widetilde h\) in \(L^2(\tau,T;L^2(\muh))\). Passing to the time-integrated form of \eqref{eq:smoothing-galerkin} shows that \(\widetilde h\) solves \eqref{eq:ratio-eq} with coefficient \(b^\lambda\) and initial datum \(\hh_0\). If \(z\) is the difference of two such solutions, then
\[
 \frac12\frac{\dd}{\dd t}\|z\|_\lambda^2+\|\nabla z\|_\lambda^2
 \le B\|z\|_\lambda\|\nabla z\|_\lambda,
\]
so uniqueness follows from Gr\"onwall's inequality. Thus \(\widetilde h=\hh\). Finally, \eqref{eq:weighted-smoothing-final} implies part~(i) after restricting to \(t\ge\tau\).
\end{proof}

% ==================================================
\section{Identification of the normal trace}
\label{app:normal-trace}
% ==================================================

We complete the proof of Lemma~\ref{lem:transmission}. Fix a time at which the conclusions established there hold. For \(s<0\), direct comparison with the interior flux gives
\begin{equation}\label{eq:parallel-normal-component}
 F(y,s)=m(y,s)
 (J^\lambda_K\cdot\vec n_s)(\Psi(y,s))
 \quad\text{for a.e. }(y,s),
\end{equation}
where \(\vec n_s(\Psi(y,s))=\vec n(y)\). Since \(m(\cdot,s)\to1\) uniformly as \(s\uparrow0\), the \(L^2(\Gamma)\)-continuity of \(F\) proves the parallel-surface convergence in \eqref{eq:g-def}.

It remains to verify the Green identity. On \(K\), one has \(w_\lambda=1\) and \(\nabla U_\lambda=0\); hence \(\rho^\lambda=\hh\), \(J_K^\lambda=J^\lambda\), and the fixed-time weak formulation gives
\[
 \divg J_K^\lambda(t)=\partial_t\rho^\lambda(t)
 =\partial_t\hh(t)\in L^2(K)
 \quad\text{in }\mathcal D'(K).
\]
Thus \(J_K^\lambda(t)\in H(\divg;K)\), and its normal trace is a priori an element of \(H^{-1/2}(\Gamma)\). We show that it is represented by \(g^\lambda(t)=F(t,\cdot,0)\in L^2(\Gamma)\).

Let \(\varphi\in C^1(\overline K)\) and fix \(\eps\in(0,\delta')\). For all sufficiently large \(k\), define
\[
 \chi_{\eps,k}(s)=
 \begin{cases}
  1,&s\le-\eps-k^{-1},\\
  k(-\eps-s),&-\eps-k^{-1}<s<-\eps,\\
  0,&s\ge-\eps.
 \end{cases}
\]
Then \(-\chi_{\eps,k}'=k\one_{(-\eps-k^{-1},-\eps)}\) has integral one. Using \(\varphi\chi_{\eps,k}(d_s)\in H^1_0(K)\) in the distributional identity for \(\divg J_K^\lambda\), and then applying the coarea formula and \eqref{eq:parallel-normal-component}, gives
\begin{align}
 &\int_K\chi_{\eps,k}(d_s)
 \bigl(J^\lambda_K\cdot\nabla\varphi
 +(\divg J^\lambda_K)\varphi\bigr)\,\dd x \notag\\
 &\qquad=-\int_{-\delta'}^0\chi_{\eps,k}'(s)
 \int_\Gamma F(y,s)\varphi(\Psi(y,s))\,\dd\Hm(y)\,\dd s .
 \label{eq:cutoff-green}
\end{align}
The surface Jacobian cancels exactly because of the definition of \(F\). The map
\[
 s\longmapsto\int_\Gamma
 F(y,s)\varphi(\Psi(y,s))\,\dd\Hm(y)
\]
is continuous. Letting \(k\to\infty\) in \eqref{eq:cutoff-green} therefore yields
\[
 \int_{\{d_s<-\eps\}}
 \bigl(J_K^\lambda\cdot\nabla\varphi
 +(\divg J_K^\lambda)\varphi\bigr)\,\dd x
 =\int_\Gamma F(y,-\eps)\varphi(\Psi(y,-\eps))\,\dd\Hm(y).
\]
Finally, \(F(\cdot,-\eps)\to g^\lambda\) in \(L^2(\Gamma)\) as \(\eps\downarrow0\), while \(\varphi(\Psi(\cdot,-\eps))\to\varphi|_\Gamma\) uniformly. Dominated convergence gives \eqref{eq:trace-green} for \(\varphi\in C^1(\overline K)\). Density and the trace theorem extend the identity to every \(\varphi\in H^1(K)\), proving that the \(H(\divg;K)\) normal trace is represented by \(g^\lambda\).

% ==================================================
\section{Convex geometry and Gaussian layer estimates}\label{app:geometry}
% ==================================================

% We collect the geometric preliminaries used in the body; the proofs are classical,
% and we indicate the argument.

% \begin{lemma}\label{lem:convex-geom}
% The domain $K$ has the following properties:
% (i) $\proj_K$ is single-valued and $1$-Lipschitz on $\Rd$, and
% $\dist^2\in C^{1,1}(\Rd)$ with $\nabla\dist^2=2(x-\proj_Kx)$;
% (ii) $\dist\in C^2(K^c)$ and
% $\nabla\dist(x)=\vec n(\proj_Kx)$ for every $x\in K^c$;
% (iii) if $x=y+s\vec n(y)$ with $y\in\partial K$ and $s>0$, then the
% tangential eigenvalues of $D^2\dist(x)$ are
% $\kappa_i(y)/(1+s\kappa_i(y))$ and its normal eigenvalue is zero.  In
% particular, $0\le\Delta\dist\le C_K$ and $|D^2\dist|\le C_K$ on all of
% $K^c$;
% (iv) if $\partial K\in C^3$, the collar map
% $\Psi(y,s)=y+s\vec n(y)$ is a $C^2$ diffeomorphism onto
% $\{|d_s|<\delta_0\}$ with Jacobian $m\in C^1$, $m(y,0)=1$,
% $c_K\le m\le C_K$.
% \end{lemma}

% \begin{proof}[Sketch]
% (i) is the standard convex projection theorem.  For a convex set with
% $C^2$ boundary, the exterior normal map
% $(y,s)\mapsto y+s\vec n(y)$ is a global $C^1$ diffeomorphism from
% $\partial K\times(0,\infty)$ onto $K^c$: uniqueness follows from the convex
% projection theorem, and its tangential differential has eigenvalues
% $1+s\kappa_i(y)>0$.  This proves (ii), and differentiating the normal map
% gives the curvature formula in (iii).  Finally, (iv) is the two-sided
% tubular-neighborhood theorem for a $C^3$ hypersurface.
% \end{proof}

\begin{lemma}\label{lem:layer}
There is $C=C(K,d,\beta)$ such that for $\lambda\in(0,1]$ and $q\ge0$,
\begin{gather*}
  \int_{K^c}\dist(x)^q\,
  e^{-\dist(x)^2/2\lambda+\beta\dist(x)}\dd x
  \ \le\ C_q\,\lambda^{(q+1)/2},\\
  \int_{K^c\cap\{\dist\ge\delta\}}
  e^{-\dist^2/2\lambda+\beta\dist}\dd x\le Ce^{-\delta^2/4\lambda}
  \ \ (0<\delta\le1).
\end{gather*}
In particular $\muh(K^c)\le C\sqrt\lambda$.
\end{lemma}

\begin{proof}
By the coarea formula with $|\nabla\dist|=1$ and the Steiner formula, $\Hm(\{\dist=s\})\le C_K(1+s)^{d-1}$, so the integral is bounded by $C\int_0^\infty s^qe^{-s^2/2\lambda+\beta s}(1+s)^{d-1}\dd s$. Completing the square, $e^{-s^2/2\lambda+\beta s}\le e^{\beta^2\lambda/2}e^{-(s-\beta\lambda)^2/2\lambda}$, and substituting $s=\beta\lambda+\sqrt\lambda\,\sigma$ gives the claim, the shift $\beta\lambda\le\beta\sqrt\lambda$ contributing only constants. For the tail, if $s\ge\delta>0$, Young's inequality gives
\[
  -\frac{s^2}{2\lambda}+\beta s
  \le-\frac{\delta^2}{4\lambda}-\frac{s^2}{8\lambda}+2\beta^2\lambda.
\]
Since $\int_0^\infty e^{-s^2/(8\lambda)}(1+s)^{d-1}\dd s\le C\sqrt\lambda$ for $0<\lambda\le1$, the stated tail bound follows with a constant independent of $\delta$.
\end{proof}

% ==================================================
\section{Difference quotients in tangential direction: proof of
Lemma~\ref{lem:tangential}}\label{app:dq}
% ==================================================

\begin{proof}[Proof of Lemma~\ref{lem:tangential}]
We give a self-contained proof. Fix a time $t$ for which $\partial_t\hh(t)\in L^2(\muh)$, $E(t)<\infty$, and the fixed-time weak formulation holds; these conditions are satisfied for almost every $t\in(0,T)$. We omit $t$ from the notation where possible and write $f:=\partial_tq(t)$, which satisfies $\int\!\!\int a f^2\dd\Hm\dd s\le C\|\partial_t\hh(t)\|_\lambda^2$ by comparability of $\muh$ with $a\dd\Hm\dd s$ on the collar.

\emph{Step 1: the fixed-time elliptic identity.} Since $\partial_t\hh\in L^2_{\mathrm{loc}}((0,T);L^2(\muh))$ by Lemma~\ref{lem:smoothing} and \eqref{eq:collar-eq} holds against test functions in $C_c^\infty((0,T)\times(\text{collar}))$, for each fixed $\phi\in C_c^\infty(\text{collar})$ the identity
\begin{equation}\label{eq:elliptic-identity}
  \int\!\!\int a\bigl(\partial_sq+\beta_sq\bigr)\partial_s\phi
  +\int\!\!\int a\bigl(\mathsf B\nG q+\beta_\Gamma q\bigr)\cdot\nG\phi
  =-\int\!\!\int a\,f\,\phi
\end{equation}
holds for a.e.\ $t$ (Lebesgue points of the $L^1_{\mathrm{loc}}$ functions of $t$ appearing on both sides); taking a countable dense family of $\phi$ and intersecting the exceptional sets, \eqref{eq:elliptic-identity} holds for a.e.\ $t$ and all $\phi\in H^1(a)$ with compact support in the collar (by density, both sides being continuous in $\phi$ for the $H^1(a)$ topology). We chose the fixed time above outside this exceptional set.

\emph{Step 2: chart notation and difference quotients.} Fix a finite atlas of $\Gamma$. In a chart $X:U\subset\R^{d-1}\to\Gamma$, let $j_X(y):=\sqrt{\det(DX(y)^\top DX(y))}$ and absorb the surface Jacobian into the weight:
\[
  \widetilde m(y,s):=m(X(y),s)j_X(y),\qquad
  a(y,s):=e^{-s_+^2/2\lambda}\widetilde m(y,s).
\]
All integrals in the chart calculation below are with respect to $\dd y\dd s$, and tangential derivatives are coordinate derivatives $\nabla_y$. We retain $q,f,\beta_s,\beta_\Gamma,\mathsf B$ for their coordinate versions, with the metric coefficients absorbed into $\mathsf B$ and $\beta_\Gamma$. The coefficients satisfy the corresponding uniform bounds in \eqref{eq:collar-coeff-bounds}, and $\mathsf B$ remains uniformly elliptic. Moreover, $c\le\widetilde m\le C$ and $\|\widetilde m\|_{C^1}\le C$ on the chart supports. The pulled-back quadratic form has no mixed normal--tangential terms because $\partial_s\Psi$ is orthogonal to the parallel surfaces. For $e\in\R^{d-1}$, $|e|=1$, and $0<|\eta|<\eta_0$ small, define
\[
  \tau_\eta v(y,s):=v(y+\eta e,s),\qquad
  D^\eta v:=\eta^{-1}(\tau_\eta v-v).
\]
We use three elementary facts, all uniform in $\lambda$ and $\eta$:
\begin{enumerate}[label=(\alph*),leftmargin=2.2em]
\item\emph{(Weight comparability and difference quotient.)} $\tau_\eta a/a=\tau_\eta\widetilde m/\widetilde m\in[c,C]$, since the singular factor $e^{-s_+^2/2\lambda}$ does not depend on $y$ and cancels in the ratio; likewise
\[
  D^\eta a=e^{-s_+^2/2\lambda}\,D^\eta\widetilde m,
  \qquad
  |D^\eta a|\le\|\nabla_y\widetilde m\|_\infty\,e^{-s_+^2/2\lambda}\le C\,a .
\]
The singular factor is \emph{never} differentiated: tangential translations commute with $s$.
\item\emph{(Discrete calculus.)} Discrete integration by parts, $\int u\,D^\eta v=-\int(D^{-\eta}u)\,v$ for $u$ or $v$ compactly supported in the chart, and the discrete Leibniz rule $D^\eta(uv)=(\tau_\eta u)D^\eta v+v\,D^\eta u$.
\item\emph{(Weighted difference-quotient bounds.)} For $v\in H^1_{\mathrm{loc}}$ with $\int a(|v|^2+|\nabla_yv|^2)<\infty$ on a neighborhood of the support of the cutoff below,
\[
  \|D^\eta v\|_{L^2(a;\,\mathrm{supp})}\le C\|\nabla_yv\|_{L^2(a)},
  \qquad
  \|D^{-\eta}w\|_{L^2(a;\,\mathrm{supp})}\le C\|\nabla_yw\|_{L^2(a)},
\]
by the fundamental theorem of calculus along the segment $[y,y+\eta e]$ and the comparability (a): the weight along the segment is $a(y+r\eta e,s)\asymp a(y,s)$ uniformly.
\end{enumerate}

\emph{Step 3: the test function.} Let $\chi=\chi(y,s)\in C_c^\infty(U\times(-\delta,\delta))$, $0\le\chi\le1$, $\chi\equiv1$ on the half-sized chart times $(-\delta',\delta')$, $|\nabla\chi|\le C$. Set $v:=D^\eta q$ and take
\[
  \phi:=-D^{-\eta}\bigl(\chi^2v\bigr)
\]
in the coordinate form of \eqref{eq:elliptic-identity}, admissible by (c) for small $\eta$ ($\phi\in H^1(a)$ with compact support, since $q\in H^1(a)$ locally by $E(t)<\infty$). By discrete integration by parts (b), and since $D^\eta$ commutes with $\partial_s$ and $\nabla_y$,
\begin{equation}\label{eq:dq-tested}
  \int\!\!\int D^\eta\bigl[a(\partial_sq+\beta_sq)\bigr]\,
  \partial_s(\chi^2v)
  +\int\!\!\int D^\eta\bigl[a(\mathsf B\nabla_y q+\beta_\Gamma q)\bigr]
  \cdot\nabla_y(\chi^2v)
  =\int\!\!\int a\,f\,D^{-\eta}(\chi^2v).
\end{equation}

\emph{Step 4: expansion and coercivity.} By the Leibniz rule (b),
\[
  D^\eta\bigl[a\mathsf B_{ij}\partial_jq\bigr]
  =(\tau_\eta a)(\tau_\eta\mathsf B_{ij})\,\partial_jv
  +\bigl[D^\eta(a\mathsf B_{ij})\bigr]\partial_jq,
  \qquad
  \bigl|D^\eta(a\mathsf B_{ij})\bigr|
  \le C\bigl(|\tau_\eta a|+|D^\eta a|\bigr)\le C\,a,
\]
using (a) and $\mathsf B\in C^1$; identically for the terms with $\beta_\Gamma q$, $\partial_sq$, $\beta_sq$ (producing $|D^\eta\beta|\le C$ by \eqref{eq:collar-coeff-bounds}). Expanding also $\nabla_y(\chi^2v)=\chi^2\nabla_y v+2\chi\nabla_y\chi\,v$ and $\partial_s(\chi^2v)=\chi^2\partial_sv+2\chi\partial_s\chi\,v$, the left-hand side of \eqref{eq:dq-tested} equals
\[
  \underbrace{\int\!\!\int(\tau_\eta a)\Bigl[
  (\tau_\eta\mathsf B)\nabla_y v\cdot\nabla_y v+|\partial_sv|^2\Bigr]\chi^2}_{
  \ge\,c_0\int\!\!\int a\,\chi^2(|\nabla_y v|^2+|\partial_sv|^2)\ \text{by (a),
  ellipticity}}
  \;+\;\mathcal E,
\]
where every term of $\mathcal E$ is of one of the forms
\[
  \int\!\!\int a\,O\bigl(|\nabla_{y,s}q|+|q|+|v|\bigr)\,
  \bigl(\chi^2|\nabla_{y,s}v|+\chi|\nabla\chi||v|\bigr),
\]
with bounded coefficients: the coefficient difference quotients contribute $a\,O(|\nabla q|+|q|)$ by the displayed bound, the drift terms contribute $a\,O(|v|+|q|)$, and no term ever contains a derivative of the singular factor. By Young's inequality and (c) (which gives $\int a\chi^2v^2\le C\int a|\nabla_y q|^2\le CE(t)$ on the enlarged support),
\[
  |\mathcal E|\le\frac{c_0}4\int\!\!\int a\,\chi^2
  \bigl(|\nabla_y v|^2+|\partial_sv|^2\bigr)
  +C\bigl(E(t)+\|\hh(t)\|_\lambda^2\bigr).
\]

\emph{Step 5: the data term.} The right side of \eqref{eq:dq-tested} is estimated by Cauchy--Schwarz and (c):
\[
  \Bigl|\int\!\!\int a f\,D^{-\eta}(\chi^2v)\Bigr|
  \le C\|f\|_{L^2(a)}\,\bigl\|\nabla_y(\chi^2v)\bigr\|_{L^2(a)}
  \le C\|f\|_{L^2(a)}\Bigl(\|\chi\nabla_y v\|_{L^2(a)}+\|v\|_{L^2(a)}\Bigr),
\]
and Young's inequality absorbs $\frac{c_0}4\|\chi\nabla_y v\|^2$, leaving $C(\|f\|_{L^2(a)}^2+E(t))$. Here $f=\partial_tq(t)$ is prescribed by the fixed-time elliptic equation and is controlled by Lemma~\ref{lem:smoothing}. The proof of that lemma, obtained by testing the ratio equation with $t\partial_t\hh$, does not use tangential regularity, so the present estimate does not depend on the conclusion of Lemma~\ref{lem:tangential}.

\emph{Step 6: conclusion.} Combining Steps 4--5 and absorbing the indicated terms,
\[
  \int\!\!\int a\,\chi^2\bigl(|\nabla_y D^\eta q|^2
  +|\partial_sD^\eta q|^2\bigr)
  \ \le\ C\Bigl(\|\partial_t\hh(t)\|_\lambda^2+E(t)
  +\|\hh(t)\|_\lambda^2\Bigr),
\]
uniformly in $\eta$ and $\lambda$. Letting $\eta\to0$, weak compactness in $L^2(a)$ and lower semicontinuity give $\nabla_y\partial_eq,\ \partial_s\partial_eq\in L^2(a)$ on the half chart times $(-\delta',\delta')$ with the same bound. Summing over the directions $e$ of a basis and over the finite atlas yields the claim for both $\nG^2q$ and $\partial_s\nG q$: the conversion to intrinsic derivatives adds only bounded-coefficient first-derivative terms, controlled by $E(t)$. (The interior side needs no separate argument: the estimate above is carried out on the full collar $(-\delta,\delta)$ at once, the normal variable being untouched by the difference quotients.)
\end{proof}

% ==================================================

\section*{Acknowledgment}
This research was supported by the National Research Foundation of Korea (NRF) grants No. RS-2023-00219980.

\section*{Disclosure of AI use}

The authors used a large language model (LLM) to assist with language editing and with checking proofs in preliminary drafts of this manuscript. The tool also assisted with the proof of Proposition~\ref{prop:sharpness} and the construction of the counterexample in Theorem~\ref{thm:counterexample}. For the latter, the authors supplied the underlying scenario: the no-flux problem starts at an unstable equilibrium, while the boundary layer amplifies an unstable mode, causing the penalized solution to move away from that equilibrium in finite time. The tool was then prompted to find a concrete example realizing this scenario. The authors reviewed and revised the AI-assisted material, independently verified all mathematical arguments and constructions, and take full responsibility for the entire content of the manuscript.

\bibliographystyle{abbrv}
\bibliography{ref}

\end{document}